\documentclass[11pt,a4paper]{article}

\usepackage[T1]{fontenc}
\usepackage{lmodern}
\usepackage[margin=1in]{geometry}
\usepackage{amsmath,amssymb,amsthm}
\usepackage{microtype}
\usepackage[hidelinks]{hyperref}
\hypersetup{
 pdftitle={Clique partitions and distinct color sets in list colorings},
 pdfauthor={Jun Yue}
}

\newtheorem{theorem}{Theorem}[section]
\newtheorem{corollary}[theorem]{Corollary}
\newtheorem{lemma}[theorem]{Lemma}
\newtheorem{proposition}[theorem]{Proposition}
\newtheorem{conjecture}[theorem]{Conjecture}
\DeclareMathOperator{\ch}{ch}
\DeclareMathOperator{\per}{per}
\newcommand{\E}{\mathbb E}
\newcommand{\PP}{\mathbb P}

\title{Clique partitions and distinct color sets in list colorings}
\author{Jun Yue\\[0.4em]
\normalsize School of Mathematical Sciences, Tiangong University\\
\normalsize Tianjin, China\\
\normalsize Email address:
\href{mailto:yuejun06@126.com}{\texttt{yuejun06@126.com}}}
\date{}

\begin{document}

\maketitle

\begin{abstract}
For an \(n\)-vertex graph \(G\) and an integer \(m\ge n\),
let \(f'(G,m)\) be the minimum number of distinct color sets
used by proper colorings from \(m\)-element lists that are
disjoint at distinct nonadjacent vertices.
A clique-partition assignment gives a common list to the
vertices of each part of a clique partition and uses disjoint
lists on different parts.
In this paper, we construct a graph on \(20\) vertices for
which these assignments do not attain \(f'(G,20)\), and a
fixed graph on \(25\) vertices for which they fail to attain
the minimum for every allowed list size.
For every triangle-free graph \(G\) on \(n\ge2\) vertices,
we prove
\[
 f'(G,m)=\binom{m}{2}^{\nu(G)}m^{n-2\nu(G)}
 \qquad (m\ge n),
\]
where \(\nu(G)\) is the size of a maximum matching, and
characterize all minimizing assignments in terms of maximum
matchings.
For clique blow-ups of each fixed odd cycle of length at
least five, we determine the minimum up to constant factors
uniformly over the blow-up size and all allowed list sizes,
and obtain an exponential separation from the clique-partition
value.
We also prove that a noncomplete graph \(G\) satisfies
\(f'(G,n)=n^2\) if and only if all edges of its complement
have a common endpoint.
These results disprove the general clique-partition conjecture
of Gravier, Maffray, and Mohar (\emph{Discrete Mathematics}, 2003)
and its extension by Wu and Zhang
(\emph{Applied Mathematics Letters}, 2006), while establishing
both triangle-free conjectures.
\end{abstract}



\raggedbottom

\section{Introduction}

All graphs considered are finite, simple, and nonempty.
Let \(G=(V,E)\) be a graph, put \(n=|V|\), and denote its
complement by \(\overline G\). A proper coloring of \(G\)
is a map \(\varphi:V\to\mathbb N\) such that
\(\varphi(u)\ne\varphi(v)\) whenever \(uv\in E\).
A list assignment \(L\) on \(G\) assigns a finite set
\(L(v)\subseteq\mathbb N\) of available colors to each
\(v\in V\). A proper \(L\)-coloring is a proper coloring
\(\varphi\) satisfying \(\varphi(v)\in L(v)\) for every
\(v\in V\). For a positive integer \(m\), the graph \(G\)
is \(m\)-choosable if every list assignment \(L\) with
\(|L(v)|=m\) for all \(v\in V\) admits a proper
\(L\)-coloring. The choice number \(\ch(G)\) is the least
such integer \(m\).

If \(H\) is a graph vertex-disjoint from \(G\), the join
\(G\vee H\) is obtained from their disjoint union by adding
all edges between \(V\) and \(V(H)\). For each positive
integer \(k\), let \(I_k\) denote the edgeless graph on
\(k\) vertices. Define
\(f(G)=\min\{k\ge1:\ch(G\vee I_k)>n\}\).
Dantas, Gravier, and Maffray introduced this parameter
in their study of equality in the bound
\(\ch(G)+\ch(\overline G)\le n+1\) of Erd\H{o}s, Rubin,
and Taylor~\cite{ERT}. Their equality characterization
first appeared in a 2000 research report, cited
in~\cite{GMM}, and was later published as~\cite{DGM}.
Bickle~\cite{B23} subsequently gave another proof of
the same characterization. Building on this formulation,
Gravier, Maffray, and Mohar~\cite{GMM} expressed \(f(G)\)
as an extremal count of the distinct color sets arising
from appropriate list assignments.

For an integer \(m\ge n\), let \(\mathcal L_m(G)\) be the
family of list assignments \(L\) on \(G\) satisfying
\[
 |L(v)|=m \quad\text{for every }v\in V,
 \qquad
 L(u)\cap L(v)=\varnothing
 \quad\text{for distinct }u,v\in V\text{ with }uv\notin E.
\]
We call these assignments admissible. Since \(m\ge n\),
every admissible assignment admits a proper coloring:
one may choose a previously unused color at each successive
vertex. Moreover, every proper coloring of an admissible
assignment is injective and therefore uses exactly \(n\)
colors. For \(L\in\mathcal L_m(G)\), put
\[
 \mathcal C_G(L)=
 \{\varphi(V):\varphi\text{ is a proper }L\text{-coloring of }G\},
 \qquad
 N_L=|\mathcal C_G(L)|,
\]
and define
\[
 f'(G,m)=\min_{L\in\mathcal L_m(G)}N_L.
\]
The parameter \(f'(G,m)\) counts each realized set once,
regardless of how many proper coloring maps have that image.
The minimum ranges over all admissible assignments of
\(m\)-element lists, so the pattern of overlaps between
the lists is part of the extremal problem.

Gravier, Maffray, and Mohar~\cite[Lemma~1]{GMM} proved
\(f(G)=f'(G,n)\); this identity is restated for later use in
Section~\ref{sec:preliminaries}. Their original definition
allows lists of size at least \(n\). Restricting each list
to an \(n\)-element subset preserves admissibility and
can only decrease the family of color sets. Consequently,
requiring lists of size exactly \(n\) yields the same
minimum.

\subsection{The clique-partition conjecture}

Let \(G\) be a graph of order \(n\), and let \(m\ge n\) be
an integer. For an admissible assignment \(L\), the vertices
whose lists contain any given color form a clique.
This suggests the following construction. Let \(\mathcal Q\)
be a partition of \(V(G)\) into nonempty cliques. Assign a
common \(m\)-element list to the vertices of each part,
choosing these lists to be pairwise disjoint.
The resulting assignment is admissible and has exactly
\(\prod_{Q\in\mathcal Q}\binom{m}{|Q|}\) color sets.
Define
\[
 w(G,m)=\min_{\mathcal Q}
             \prod_{Q\in\mathcal Q}\binom{m}{|Q|},
 \qquad
 w(G)=w(G,n),
\]
where the minimum is taken over all clique partitions
\(\mathcal Q\) of \(G\). As observed in~\cite[Theorem~2]{GMM}
and~\cite[Theorem~2.4]{WZ}, these assignments give
\(f'(G,m)\le w(G,m)\), and consequently \(f(G)\le w(G)\).
Theorem~\ref{thm:known-bounds} records these upper bounds
together with the square lower bound that will be used later.
Gravier, Maffray, and Mohar conjectured that equality
always holds in the latter bound.

\begin{conjecture}[Gravier, Maffray, and Mohar~{\cite[Conjecture~1]{GMM}}]
\label{conj:gmm-general}
For every graph \(G\), we have \(f(G)=w(G)\).
\end{conjecture}

Gravier, Maffray, and Mohar established this equality
for disjoint unions of cliques~\cite[Claim~2]{GMM}
and recorded \(f(C_4)=36\) and \(f(C_5)=500\)
in the discussion following~\cite[Lemma~1]{GMM}.
Wu and Zhang proposed the corresponding equality
for all list sizes \(m\ge n\).

\begin{conjecture}[Wu and Zhang~{\cite[Conjecture~2.5]{WZ}}]
\label{conj:general}
For every \(n\)-vertex graph \(G\) and every integer
\(m\ge n\), we have \(f'(G,m)=w(G,m)\).
\end{conjecture}

Our counterexamples are clique blow-ups of odd cycles.
For an odd integer \(\ell\ge5\) and an integer \(t\ge1\),
let \(G_{\ell,t}=C_\ell[K_t]\) be the graph obtained by
replacing the vertices of \(C_\ell\) with pairwise disjoint
\(t\)-vertex cliques \(V_0,\ldots,V_{\ell-1}\).
Vertices in distinct classes \(V_i\) and \(V_j\) are adjacent
if and only if \(j-i\equiv\pm1\pmod{\ell}\).

\begin{theorem}\label{thm:counterexamples}
The graph \(G_{5,4}\) on \(20\) vertices satisfies
\[
 f(G_{5,4})
 \le \binom{50}{20}
 < \binom{20}{8}^{2}\binom{20}{4}
 =w(G_{5,4}).
\]
The graph \(G_{5,5}\) on \(25\) vertices satisfies
\[
 f'(G_{5,5},m)<w(G_{5,5},m)
 \qquad\text{for every integer }m\ge25.
\]
\end{theorem}

The first assertion disproves
Conjecture~\ref{conj:gmm-general} and hence also
Conjecture~\ref{conj:general}. The second strengthens this
failure by keeping the graph fixed while allowing every
permitted list size. Thus the equality \(f(C_5)=w(C_5)\)
does not extend to all clique blow-ups of \(C_5\).
In both examples, colors are shared on the overlapping
cliques \(V_i\cup V_{i+1}\), with indices taken modulo \(\ell\).
This overlap permits a smaller family of color sets than
the families obtained by giving disjoint palettes to the
parts of any clique partition.

\subsection{The triangle-free conjecture}

Suppose now that \(G\) is triangle-free and has \(n\ge2\)
vertices. Let \(\nu(G)\) denote the size of a maximum matching
in \(G\). Every clique partition of \(G\) consists of edges
forming a matching and singleton parts corresponding to
unmatched vertices. A partition containing \(k\) edges has
value \(\binom{m}{2}^{k}m^{n-2k}\). Since
\(\binom{m}{2}<m^2\), this value is minimized when
\(k=\nu(G)\). Thus
\(w(G,m)=\binom{m}{2}^{\nu(G)}m^{n-2\nu(G)} \qquad (m\ge n\ge2)\).
Gravier, Maffray, and Mohar stated the triangle-free case
of Conjecture~\ref{conj:gmm-general} as a separate conjecture.

\begin{conjecture}[Gravier, Maffray, and Mohar~{\cite[Conjecture~2]{GMM}}]
\label{conj:gmm-trianglefree}
For every triangle-free graph \(G\) on \(n\ge2\) vertices,
\[
 f(G)=\binom{n}{2}^{\nu(G)}n^{n-2\nu(G)}.
\]
\end{conjecture}

Gravier, Maffray, and Mohar proved this conjecture for
a recursively constructed class of trees~\cite[Corollary~1]{GMM}.
Wu and Zhang subsequently extended its formulation to
all allowed list sizes and proved the corresponding
formula for forests~\cite[Theorem~3.3]{WZ}.
That known formula is restated in
Theorem~\ref{thm:forest-known}; the proposed extension
is the following.

\begin{conjecture}[Wu and Zhang~{\cite[Conjecture~2.6]{WZ}}]
\label{conj:trianglefree}
For every triangle-free graph \(G\) on \(n\ge2\) vertices
and every integer \(m\ge n\), we have \(f'(G,m)=w(G,m)\).
\end{conjecture}

Cycles introduce multiplicities that must be controlled
to pass from the forest case to general triangle-free
graphs. Although admissibility restricts each color to
at most two lists, several proper coloring maps may
still have the same image. We therefore estimate the
multiplicity of a prescribed color set together with
the probability that independent uniform choices from
the lists repeat a color. This joint estimate yields
the conjectured minimum and also determines all the
assignments that attain it.

\begin{theorem}\label{thm:trianglefree}
Let \(G\) be a triangle-free graph on \(n\ge2\) vertices,
and let \(m\ge n\) be an integer. Then
\[
 f'(G,m)=w(G,m)
       =\binom{m}{2}^{\nu(G)}m^{n-2\nu(G)}.
\]
An assignment \(L\in\mathcal L_m(G)\) attains this minimum
if and only if there is a maximum matching \(M\) in \(G\)
such that \(L(u)=L(v)\) for every \(uv\in M\), and the
common lists on the matched pairs and the lists on the
unmatched vertices are pairwise disjoint.
For the one-vertex graph,
\(f'(K_1,m)=w(K_1,m)=m\) for every integer \(m\ge1\).
\end{theorem}

The specialization \(m=n\) proves
Conjecture~\ref{conj:gmm-trianglefree}, and the assertion
for arbitrary allowed list sizes establishes
Conjecture~\ref{conj:trianglefree}. The equality statement
also specifies exactly how the lists in a minimizing
assignment can overlap: their common palettes are supported
on the edges of a maximum matching.
Together with Theorem~\ref{thm:counterexamples}, this locates
the failure of the general conjecture beyond the triangle-free
class, throughout which clique-partition assignments remain
optimal. The directed-cycle estimates and fractional matching
arguments behind the proof are outlined in
Section~\ref{sec:outline}.

\subsection{Further results and related work}

Fix an odd integer \(\ell\ge5\), and let \(t\ge1\) and
\(m\ge\ell t\) be integers. Every clique of \(G_{\ell,t}\)
has at most \(2t\) vertices, so each color in an admissible
assignment occurs in at most \(2t\) lists. Counting
vertex-color incidences therefore shows that the union of
the lists contains at least \(\lceil\ell m/2\rceil\) colors.
Our palette construction attains this bound. For even \(m\),
it uses pairwise disjoint palettes of size \(m/2\) on the
cliques \(V_i\cup V_{i+1}\); for odd \(m\), the palette sizes
are adjusted and one additional color is supported on a
single class. Since every proper coloring uses exactly
\(\ell t\) colors, this construction gives
\(f'(G_{\ell,t},m) \le \binom{\lceil\ell m/2\rceil}{\ell t}\).
The following theorem shows that this upper bound is optimal
up to a constant factor depending only on \(\ell\).

\begin{theorem}\label{thm:blowup}
For each fixed odd integer \(\ell\ge5\),
\[
 f'(G_{\ell,t},m)
 =\Theta_\ell\!\left(
       \binom{\lceil\ell m/2\rceil}{\ell t}\right),
\]
uniformly over integers \(t\ge1\) and \(m\ge\ell t\).
The implicit positive constants depend only on \(\ell\).
\end{theorem}

For these graphs, the clique-partition value is
\(w(G_{\ell,t},m) =\binom{m}{2t}^{(\ell-1)/2}\binom{m}{t}\).
Combining this formula with Theorem~\ref{thm:blowup} gives
a constant \(c_\ell>0\) such that
\[
 \frac{w(G_{\ell,t},m)}{f'(G_{\ell,t},m)}
 \ge c_\ell\,2^t t^{-(\ell-1)/4}
 \qquad (t\ge1,\ m\ge\ell t).
\]
Thus the ratio grows exponentially in \(t\), uniformly in
the list size. In particular, there is \(t_0=t_0(\ell)\)
such that \(G_{\ell,t}\) is a counterexample to
Conjecture~\ref{conj:general} for every integer \(m\ge\ell t\)
whenever \(t\ge t_0\).

We also establish a stability result for assignments whose
number of color sets is within a fixed factor of the minimum.
More precisely, fix \(K\ge1\) and \(M\ge\ell\).
If \(\ell t\le m\le Mt\) and
\(L\in\mathcal L_m(G_{\ell,t})\) satisfies
\(N_L\le Kf'(G_{\ell,t},m)\), then \(L\) can be transformed
into an assignment of the palette construction by inserting
or deleting \(O_{\ell,K,M}(t)\) vertex-color incidences.
Here an incidence edit inserts or deletes one color
in the list of one vertex.

A further extremal question concerns the smallest possible
value of \(f(G)\) for a noncomplete graph \(G\) of order \(n\).
Gravier, Maffray, and Mohar proved \(f(G)\ge n^2\), with
equality whenever all edges of \(\overline G\) have a common
endpoint~\cite[Theorem~1 and Claim~1]{GMM}.
They observed that no other equality examples were known.
Theorem~\ref{thm:star} proves the converse: a noncomplete
graph \(G\) satisfies \(f(G)=n^2\) if and only if all edges
of \(\overline G\) have a common endpoint.

Related bounds arise in the study of generalized Ohba
numbers. Cano et al.~\cite[Lemma~20]{Cano} established a
product lower bound involving pairwise disjoint lists.
Writing \(\alpha(G)\) for the independence number of \(G\),
their bound applied to the lists on a maximum independent
set gives
\(f'(G,m)\ge m^{\alpha(G)} \qquad (\alpha(G)\ge2,\ m\ge n)\).
For \(G_{\ell,t}\), we have \(\alpha(G_{\ell,t})=(\ell-1)/2\),
so this yields the lower bound \(m^{(\ell-1)/2}\).
When \(m=\ell t\), Theorem~\ref{thm:blowup} determines
the exponential rate and the polynomial prefactor up to
constant factors. The product estimate uses only the sizes of the lists in
the specified disjoint subfamily, whereas the blow-up bound accounts for the
overlaps among all the lists in the assignment.

Counting color sets also connects the problem with earlier
work on the number of bases of bicircular matroids. Let \(H\) be a connected graph
with \(N\) vertices and \(M\) edges, containing a cycle,
and let \(b(H)\) and \(\tau(H)\) denote its numbers of
bicircular bases and spanning trees, respectively.
Gim\'enez, de Mier, and Noy~\cite[Theorems~2.1 and~2.2]{GMN05}
proved \(b(H)\ge\frac{M-N+1}{N}\tau(H)\) and, when
\(H\) is simple, \(b(H)\le2.42^N\tau(H)\).
These estimates relate the two counts on each fixed graph.
In the present triangle-free problem, the color-support
multigraph varies with the list assignment. The prescribed
list sizes and the allowed supports constrain that variation,
and the objective is to minimize its number of bases.

\subsection{Proof overview and organization}\label{sec:outline}

The counterexample argument has two separate counting steps.
For \(G_{\ell,t}\), palettes supported on the overlapping
cliques \(V_i\cup V_{i+1}\) attain the smallest possible
union of lists. We then compare the number of subsets of
that union with the minimum product coming from a clique
partition. This comparison proves
Theorem~\ref{thm:counterexamples}. For the fixed
\(25\)-vertex graph, a factor-by-factor comparison of
binomial coefficients makes the argument valid uniformly
over every allowed list size.

For triangle-free graphs, the support of a color has one
or two vertices, so it can be represented by a loop or
an edge of a multigraph on \(V(G)\). The realizable images
correspond to spanning subgraphs with unicyclic components.
In each component, a nonloop cycle contributes a factor
of two to the number of coloring maps with the same image.
Independent choices from the lists give a random map whose
off-diagonal transition probabilities form a fractional
matching. We estimate its cycles and color collisions
jointly, then round the extreme points in a convex
decomposition of this fractional matching to obtain a
matching. The quantitative remainder retained in this
comparison provides the strictness needed to identify
every equality case.

For odd-cycle blow-ups, the permanent-capacity inequality
estimates the total number of injective choices. The other
part of the argument controls how many choices can have one
prescribed image: after enlarging the color supports to
consecutive pairs of clique classes, this becomes a problem
of allocating colors to the classes. Uniform binomial
probability bounds give the required multiplicity estimate.
Its combination with the capacity bound leaves a nonnegative
remainder, which subsequently controls the incidence edits
in the stability result. The square lower bound is treated
separately using the disjoint-list product estimate and
Hall's theorem to characterize equality.

The remainder of the paper is organized as follows.
Section~\ref{sec:preliminaries} collects the results on
graph joins, transversals, fractional matchings, permanents,
and probability estimates used in the proofs.
Theorems~\ref{thm:counterexamples}, \ref{thm:trianglefree},
and~\ref{thm:blowup} are proved in
Sections~\ref{sec:counterexamples}, \ref{sec:trianglefree},
and~\ref{sec:blowups}, respectively.
Section~\ref{sec:star} completes the main argument by
characterizing equality in the bound \(f(G)\ge n^2\).
The appendices provide the supporting exact counts:
Appendix~\ref{app:support-count} gives the formula in terms
of color supports, and Appendix~\ref{sec:finite-count}
evaluates the \(20\)-vertex construction exactly.

\section{Preliminaries}\label{sec:preliminaries}

The results in this section are arranged according to their
roles in the proofs. We begin with the join identity, the
known extremal bounds, and the descriptions of transversals
and fractional matchings. The permanent inequalities then
yield the counting estimates in Lemma~\ref{lem:injective}
and Proposition~\ref{prop:general}, and we include their
derivations in the notation of admissible lists.
The binomial estimates, network principle, and orientation
count supply the remaining inputs for the two main
counting arguments. All the classical results are stated
with their sources in the forms used below. All logarithms
are natural, and we use the convention \(0\log0=0\).

\subsection{Graph joins, transversals, and matchings}

The following identity expresses the join threshold \(f(G)\)
as the minimum number of color sets arising from admissible
list assignments.

\begin{lemma}[Gravier, Maffray, and Mohar~{\cite[Lemma~1]{GMM}}]
\label{lem:join}
Let \(G\) be a graph, and put \(n=|V(G)|\). Then
\(f(G)=f'(G,n) =\min_{L\in\mathcal L_n(G)}|\mathcal C_G(L)|\).
\end{lemma}

We also record the known extremal bounds and the forest
formula. The clique-partition bound provides the benchmark
for the counterexamples, while the forest formula gives
the previously settled case of the triangle-free problem.
The square lower bound and its equality examples will be
the starting point for the necessity argument in
Theorem~\ref{thm:star}.

\begin{theorem}[Clique-partition and square bounds]
\label{thm:known-bounds}
Let \(G\) be a graph of order \(n\), and let \(m\ge n\)
be an integer. Then \(f'(G,m)\le w(G,m)\)
\cite[Theorem~2.4]{WZ}; for \(m=n\), this gives
\(f(G)\le w(G)\), as proved in~\cite[Theorem~2]{GMM}.
If \(G\) is noncomplete, then \(f(G)\ge n^2\), with
equality whenever all edges of \(\overline G\) have a
common endpoint~\cite[Theorem~1 and Claim~1]{GMM}.
\end{theorem}

\begin{theorem}[Wu and Zhang~{\cite[Theorem~3.3]{WZ}}]
\label{thm:forest-known}
Let \(G\) be a forest of order \(n\ge2\), and let
\(m\ge n\) be an integer. Then
\(f'(G,m)=\binom{m}{2}^{\nu(G)}m^{n-2\nu(G)}\).
For the one-vertex graph, \(f'(G,m)=m\).
\end{theorem}

For a list assignment \(L\) on a graph \(G\), write
\(\mathcal U(L)=\bigcup_{v\in V(G)}L(v)\)
for its set of colors. The bipartite incidence graph of
\(L\) has parts \(V(G)\) and \(\mathcal U(L)\), with
\(v\) adjacent to \(c\) precisely when \(c\in L(v)\).
If \(L\) is admissible, its proper colorings correspond
to matchings in this incidence graph that saturate \(V(G)\),
or equivalently to systems of distinct representatives
of the lists. We use the following two classical results
to describe their existence and their images.

\begin{theorem}[Hall~{\cite[Theorem~1]{Hall}}]
\label{thm:hall}
Let \((A_i)_{i\in I}\) be a finite family of finite sets.
There exists an injective map \(a:I\to\bigcup_{i\in I}A_i\)
with \(a(i)\in A_i\) for every \(i\in I\) if and only if
\[
 \left|\bigcup_{i\in J}A_i\right|\ge |J|
 \qquad\text{for every }J\subseteq I.
\]
Equivalently, a finite bipartite graph with parts \(I,U\)
has a matching saturating \(I\) if and only if
\(|N(J)|\ge|J|\) for every \(J\subseteq I\).
\end{theorem}

\begin{lemma}[Edmonds and Fulkerson~{\cite[Section~1]{EF}}]
\label{lem:transversal}
Let \((A_i)_{i\in I}\) be a finite family of subsets of
a finite set \(U\). The subsets \(S\subseteq U\) that
can be matched injectively to distinct indices in \(I\),
with each element belonging to the set indexed by its
matched index, are the independent sets of a matroid
on \(U\). If the family has a system of distinct
representatives, this matroid has rank \(|I|\), and
its bases are exactly the images of such systems.
\end{lemma}

An admissible assignment has a system of distinct
representatives, so Lemma~\ref{lem:transversal} identifies
its color sets with the bases of a transversal matroid
on \(\mathcal U(L)\). The incidence graph records which vertex receives
each representative, whereas the matroid records only the
chosen colors. The permanent bounds below compare these
counts by controlling how many representative systems
can have one image.

The following product bound controls the number of distinct
images using a pairwise disjoint subfamily of lists.
It allows unequal list sizes and therefore applies after
a color is forced at one vertex and removed elsewhere.
Its two-vertex case will be used in the equality
classification in Section~\ref{sec:star}.

\begin{lemma}[Cano et al.~{\cite[Lemma~20]{Cano}}]
\label{lem:two-list-bound}
Let \(H\) be a graph on \(q\ge2\) vertices, and let \(J\)
be a list assignment satisfying \(|J(w)|\ge q\) for every
\(w\in V(H)\). Let \(a_1,\ldots,a_s\) be distinct
vertices, where \(2\le s\le q\), such that their lists
are pairwise disjoint. Then
\[
 \bigl|
 \{\varphi(V(H)):
       \varphi\text{ is a proper }J\text{-coloring of }H\}
 \bigr|
 \ge\prod_{i=1}^s |J(a_i)|.
\]
In particular, two nonadjacent vertices \(a,b\) with
disjoint lists give the lower bound \(|J(a)|\,|J(b)|\).
The statement above is Lemma~20 of the third arXiv version
of~\cite{Cano}, dated 23 February 2026; it does not require
admissibility.
\end{lemma}

The overlap sizes in the triangle-free argument will be
treated as fractional matchings. Let \(G=(V,E)\) be a graph.
Its fractional matching polytope is
\[
 \mathcal P_f(G)=
 \left\{
 x\in\mathbb R_{\ge0}^{E}:
 \sum_{\substack{e\in E\\e\ni v}}x_e\le1
 \text{ for every }v\in V
 \right\}.
\]
The elements of \(\mathcal P_f(G)\) are called fractional
matchings. We use the following standard description of
its extreme points.

\begin{lemma}[Fractional matching extreme points; {\cite[Chapter~7]{LP}}]\label{lem:fractional}
Every extreme point \(x\) of \(\mathcal P_f(G)\) has
coordinates in \(\{0,1/2,1\}\). The edges with \(x_e=1\)
form a matching, and the edges with \(x_e=1/2\) form
vertex-disjoint odd cycles, disjoint from the vertices
of that matching. If \(G\) is bipartite, every extreme
point is the indicator vector of a matching.
\end{lemma}

Since \(\mathcal P_f(G)\) is a polytope, its elements
admit convex decompositions into extreme points.
We use this fact in probabilistic form: for each
\(x\in\mathcal P_f(G)\), we may choose a random extreme
point \(X\) with \(\mathbb E X=x\).
When \(G\) is bipartite, \(X\) can therefore be regarded
as the indicator vector of a randomly chosen matching.

\subsection{Permanent estimates and color-set bounds}

The next estimates count the injective choices themselves
and bound the number that can realize one prescribed image.
Both tasks can be expressed in terms of permanents of
incidence matrices.
For a positive integer \(q\) and a nonnegative
\(q\times q\) matrix \(A=(A_{ij})\), its permanent is
\(\per A=\sum_{\sigma\in S_q}\prod_{i=1}^q A_{i,\sigma(i)}\),
where \(S_q\) is the set of permutations of
\(\{1,\ldots,q\}\). Define its capacity by
\[
 \operatorname{cap}(A)=
 \inf_{z_1,\ldots,z_q>0}
 \frac{\displaystyle
       \prod_{i=1}^q
       \left(\sum_{j=1}^q A_{ij}z_j\right)}
      {\displaystyle\prod_{j=1}^q z_j}.
\]
Completing the full incidence matrix to a square matrix
permits the lower bound to count injective choices.
For a prescribed image, transposing its incidence matrix
makes the color-occurrence numbers available as row sums
in the upper bound.

\begin{theorem}[Gurvits~{\cite[Corollary~2.5 and Example~2.6]{Gurvits}}]
\label{thm:capacity}
For every nonnegative \(q\times q\) matrix \(A\),
\(\per A\ge(q!/q^q)\operatorname{cap}(A)\).
\end{theorem}

\begin{theorem}[Br\'egman~{\cite{Bregman}}]
\label{thm:bregman}
Let \(A\) be a \(q\times q\) matrix with entries in
\(\{0,1\}\) and positive row sums \(r_1,\ldots,r_q\).
Then \(\per A\le\prod_{i=1}^q(r_i!)^{1/r_i}\).
If a row sum is zero, then \(\per A=0\).
\end{theorem}

Radhakrishnan~\cite{R97} gives an entropy proof of
Theorem~\ref{thm:bregman}. Sra~\cite[Theorem~2.5 and
Remark~2.8]{Sra} discusses the equivalent permutation-entropy
formulation of Theorem~\ref{thm:capacity}. We use the
capacity formulation for the list incidence matrix.
Completing that matrix and choosing suitable doubly stochastic
weights gives the following lower bound on the total number
of proper colorings.

\begin{lemma}[A consequence of Gurvits's inequality]
\label{lem:injective}
Let \(G\) be a graph of order \(n\), let \(m\ge n\) be
an integer, and let \(L\in\mathcal L_m(G)\).
For each \(c\in\mathcal U(L)\), put
\(d_c=|\{v\in V(G):c\in L(v)\}|, \qquad p_c=\frac{d_c}{m}\).
If \(T\) is the number of proper \(L\)-colorings, then
\[
 \log T\ge
 n\log m-n
 -\sum_{c\in\mathcal U(L)}(1-p_c)\log(1-p_c).
\]
\end{lemma}

\begin{proof}
Write \(\mathcal U=\mathcal U(L)\) and \(q=|\mathcal U|\).
Since each list has size \(m\), we have \(q\ge m\ge n\).
Counting vertex-color incidences gives
\(0<p_c\le1 \quad(c\in\mathcal U), \qquad \sum_{c\in\mathcal U}p_c=n\).

Suppose first that \(q=n\). Then \(m=n\), and every list
equals \(\mathcal U\). Every proper coloring is injective,
so \(T=n!\). Since
\(\log(n!)\ge\int_1^n\log x\,dx =n\log n-n+1\),
the required inequality follows.

We may therefore assume \(q>n\), and put \(b=q-n\ge1\).
Append \(b\) all-one rows to the \(n\times q\) incidence
matrix of \(L\), obtaining a \(q\times q\) matrix \(A\).
Each proper \(L\)-coloring selects distinct columns for
the original \(n\) rows. The remaining \(b\) columns can
be assigned bijectively to the appended rows in \(b!\)
ways. Conversely, every permutation contributing to
\(\per A\) restricts to such a coloring. Hence
\(\per A=b!T\).

Index the columns by \(\mathcal U\), and define a
\(q\times q\) matrix \(X\) by
\[
 X_{ic}=
 \begin{cases}
 A_{ic}/m, & 1\le i\le n,\\[2mm]
 (1-p_c)/b, & n<i\le q.
 \end{cases}
\]
Each original row sums to one. Since
\(\sum_{c\in\mathcal U}(1-p_c)=q-n=b\), each appended
row also sums to one. The sum of column \(c\) is
\(p_c+(1-p_c)=1\). Thus \(X\) is doubly stochastic.

Whenever \(X_{ic}>0\), we have \(A_{ic}=1\).
Consequently, for positive \(z_c\), the weighted
arithmetic--geometric mean inequality gives
\[
 \sum_{c\in\mathcal U}A_{ic}z_c
 \ge
 \prod_{\substack{c\in\mathcal U\\X_{ic}>0}}
       \left(\frac{z_c}{X_{ic}}\right)^{X_{ic}}.
\]
Multiplying over the rows and using the column sums of
\(X\), we obtain
\[
 \prod_{i=1}^q
       \left(\sum_{c\in\mathcal U}A_{ic}z_c\right)
 \ge
 \exp\!\left(
       -\sum_{i=1}^q\sum_{c\in\mathcal U}
          X_{ic}\log X_{ic}\right)
 \prod_{c\in\mathcal U}z_c.
\]
The definition of \(X\) gives
\[
 -\sum_{i=1}^q\sum_{c\in\mathcal U}X_{ic}\log X_{ic}
 =
 n\log m+b\log b
 -\sum_{c\in\mathcal U}(1-p_c)\log(1-p_c).
\]
It follows that
\[
 \operatorname{cap}(A)
 \ge
 m^n b^b
 \exp\!\left(
       -\sum_{c\in\mathcal U}(1-p_c)\log(1-p_c)\right).
\]
Applying Theorem~\ref{thm:capacity} and using \(\per A=b!T\),
we conclude that
\[
 T\ge
 \frac{q!\,b^b}{q^q b!}\,m^n
 \exp\!\left(
       -\sum_{c\in\mathcal U}(1-p_c)\log(1-p_c)\right).
\]

Finally, the inequality
\(\log(1-x)\ge-x/(1-x)\), valid for \(0\le x<1\), yields
\[
 \log\frac{q!\,b^b}{q^q b!}
 =
 \sum_{j=b+1}^{q}(j-1)\log\left(1-\frac1j\right)
 \ge -(q-b)=-n.
\]
Taking logarithms in the preceding lower bound for \(T\)
proves the assertion.
\end{proof}

To express the counting bound in entropy notation, for
\(x\in[0,1]\) define
\(h(x)=-x\log x-(1-x)\log(1-x)\).
For \(0<z<1\), define
\(g(z)=-\frac{(1-z)\log(1-z)}{z}\),
and extend \(g\) continuously to \(z=1\) by \(g(1)=0\).
For \(0<z<1\), we have
\(g'(z)=\frac{\log(1-z)+z}{z^2}<0\).
Thus \(g\) is decreasing on \((0,1]\).

We apply Lemma~\ref{lem:injective} to the full list assignment
and Theorem~\ref{thm:bregman} to the incidence matrix of a
fixed image. Comparing the resulting estimates gives a lower
bound for \(f'(G,m)\) in terms of the clique number of \(G\).
The derivation also retains a nonnegative remainder that
depends on the individual color supports and will be used
in the blow-up argument.
\begin{proposition}[A consequence of the permanent estimates]\label{prop:general}
Let \(G\) be a graph, and put \(n=|V(G)|\) and
\(r=\omega(G)\). For an integer \(m\ge n\), set \(p=r/m\).
Then
\[
 \log f'(G,m)\ge
 \frac{n}{p}h(p)
 -n\left(1+\frac{\log(r!)}{r}-\log r\right).
\]
\end{proposition}

\begin{proof}
Fix \(L\in\mathcal L_m(G)\), and write
\(\mathcal U=\bigcup_{v\in V(G)}L(v)\).
For each \(c\in\mathcal U\), let
\(d_c=|\{v\in V(G):c\in L(v)\}|, \qquad p_c=\frac{d_c}{m}\).
Admissibility implies that the vertices whose lists contain
\(c\) form a clique. Consequently,
\(1\le d_c\le r, \qquad 0<p_c\le p\le1, \qquad \sum_{c\in\mathcal U}p_c=n\).

Let \(T\) be the total number of proper \(L\)-colorings.
For \(S\in\mathcal C_G(L)\), let \(T(S)\) be the number
of these colorings with image \(S\).
Every proper coloring is injective, so \(|S|=n\).
Let \(A_S\) be the \(n\times n\) matrix whose rows are
indexed by \(V(G)\), whose columns are indexed by \(S\),
and whose entries satisfy
\[
 (A_S)_{vc}=
 \begin{cases}
 1, & c\in L(v),\\
 0, & c\notin L(v).
 \end{cases}
\]
The bijections from \(V(G)\) to \(S\) counted by
\(\per A_S\) are precisely the proper \(L\)-colorings
with image \(S\). Thus \(T(S)=\per A_S\).

The column of \(A_S\) indexed by \(c\) has sum \(d_c\).
Applying Theorem~\ref{thm:bregman} to \(A_S^{\mathsf T}\)
therefore gives
\[
 T(S)=\per A_S^{\mathsf T}
 \le \prod_{c\in S}(d_c!)^{1/d_c};
\]
the theorem applies because these column sums are positive.
The sequence \(\log(d!)/d\) is nondecreasing for positive
integers \(d\), since it is the mean of the first \(d\)
terms of the increasing sequence
\(\log1,\log2,\ldots\).
Using \(d_c\le r\) and \(|S|=n\), we obtain
\(T(S)\le(r!)^{n/r}\).
Summing over all realized color sets yields
\(T=\sum_{S\in\mathcal C_G(L)}T(S) \le N_L(r!)^{n/r}\).
Hence
\(\log N_L\ge\log T-\frac{n}{r}\log(r!)\).

The monotonicity of \(g\) gives
\(g(p_c)\ge g(p)\) for every \(c\in\mathcal U\).
Lemma~\ref{lem:injective} gives
\(\log T \ge n\log m-n +\sum_{c\in\mathcal U}p_cg(p_c)\).
Combining this with the multiplicity estimate, and using
\(m=r/p\), \(\sum_c p_c=n\), and
\(\frac{h(p)}{p}=-\log p+g(p)\),
we obtain
\[
 \begin{aligned}
 \log N_L
 &\ge n\log m-n
      +\sum_{c\in\mathcal U}p_cg(p_c)
      -\frac{n}{r}\log(r!)\\
 &=\frac{n}{p}h(p)
   -n\left(1+\frac{\log(r!)}{r}-\log r\right)\\
 &\qquad
   +\sum_{c\in\mathcal U}
       p_c\bigl(g(p_c)-g(p)\bigr).
 \end{aligned}
\]
The last sum is nonnegative. The resulting lower bound
is independent of \(L\), so taking the minimum over
\(L\in\mathcal L_m(G)\) proves the proposition.
\end{proof}

\subsection{Probability estimates and orientations}

The odd-cycle blow-up argument requires estimates that remain
uniform as both the list size and the class size vary.
We use the following consequences of Stirling's formula
and the binomial concentration bound. The constants in
the Stirling estimates are absolute.

\begin{lemma}[Stirling estimates; Robbins~\cite{Robbins}]
\label{lem:stirling}
For every integer \(j\ge1\),
\[
 j!=\sqrt{2\pi j}\,(j/e)^j e^{\rho_j},
 \qquad
 \frac{1}{12j+1}<\rho_j<\frac{1}{12j}.
\]
Consequently, there is an absolute constant \(C\) such that
\(j!\le C\sqrt j\,(j/e)^j\) for every \(j\ge1\), and
\(\max_a\PP(\operatorname{Bin}(j,1/2)=a)\le Cj^{-1/2}\).
For integers \(1\le k<N\), another consequence is the
uniform estimate
\[
 \log\binom Nk
 =Nh(k/N)-\frac12\log\bigl(k(1-k/N)\bigr)+O(1),
\]
where the constant in \(O(1)\) is independent of \(N,k\).
\end{lemma}

\begin{lemma}[A binomial form of Hoeffding's inequality~{\cite[Theorem~2]{Hoeffding}}]
\label{lem:hoeffding}
Let \(Z\sim\operatorname{Bin}(N,1/2)\), where \(N\ge1\).
For every \(a>0\),
\(\PP(|Z-N/2|\ge a)\le2\exp(-2a^2/N)\).
\end{lemma}

For the directed-cycle estimate, we will also compare
an escape probability with the energy of a finite network.
We use the following variational principle, with each
undirected edge contributing once to the energy.

\begin{lemma}[Dirichlet principle~{\cite[Chapter~2, Exercise~2.13]{LyP}}]
\label{lem:dirichlet}
Let \(V\) be a finite vertex set, with symmetric conductances
\(c_{xy}=c_{yx}\ge0\) and \(c_{xx}=0\). Fix distinct vertices
\(u,v\), and put
\(\mathcal E(g)=\frac12\sum_{x,y\in V}c_{xy}(g(x)-g(y))^2\).
If \(h(u)=0\), \(h(v)=1\), and
\(\sum_y c_{xy}(h(x)-h(y))=0\) for every
\(x\notin\{u,v\}\), then
\[
 \mathcal E(h)=
 \min\{\mathcal E(g):g(u)=0,\ g(v)=1\}.
\]
On any positive-conductance component meeting neither
boundary vertex, one may choose \(h\) to be constant.
\end{lemma}

An orientation with one outgoing edge at each vertex
relates representative systems to unicyclic components.
We use the following standard simple-graph form of this
correspondence, as described in~\cite[Section~2]{GJ}.

\begin{lemma}[Cycle orientations; {\cite[Section~2]{GJ}}]
\label{lem:cycle-orientations}
Let \(H\) be a finite simple graph. It has an orientation
with outdegree one at every vertex if and only if every
connected component is unicyclic. If there are \(k\)
components, there are exactly \(2^k\) such orientations:
the cycle in each component has two cyclic orientations,
and every edge off the cycle points towards it.
\end{lemma}

In the list model of Section~\ref{sec:trianglefree}, individual
colors label the edges, and self-loops and parallel edges
are allowed. The correspondence is verified there with
these conventions: the distinction between two labeled
parallel edges and two selections of the same edge is
essential when multiplicities and collisions are counted.

\section{Counterexamples}\label{sec:counterexamples}

Fix an integer \(k\ge2\), and put \(\ell=2k+1\).
For an integer \(t\ge1\), let \(V_0,\ldots,V_{\ell-1}\)
be the \(t\)-vertex clique classes of \(G_{\ell,t}\)
defined in the introduction. Class indices are taken
modulo \(\ell\).
We begin by determining two quantities for the same graph:
the smallest union of lists attainable by the palette
construction and the exact clique-partition value. Their
comparison will give both counterexamples.

\begin{lemma}\label{lem:palettes}
For every pair of integers \(t\ge1\) and \(m\ge\ell t\),
\[
 f'(G_{\ell,t},m)
 \le \binom{\lceil\ell m/2\rceil}{\ell t},
 \qquad
 w(G_{\ell,t},m)
 =\binom{m}{2t}^{k}\binom{m}{t}.
\]
Moreover, the minimum possible cardinality of the union
of the lists is
\[
 \min_{L\in\mathcal L_m(G_{\ell,t})}
 \left|\bigcup_{v\in V(G_{\ell,t})}L(v)\right|
 =\left\lceil\frac{\ell m}{2}\right\rceil.
\]
\end{lemma}

\begin{proof}
We first determine the minimum possible size of the union
of the lists. The classes met by a clique in \(G_{\ell,t}\)
correspond to a clique in \(C_\ell\). Since \(\ell\ge5\),
every clique therefore meets at most two classes, which
must be consecutive if both are met. In particular,
every clique has at most \(2t\) vertices.
Let \(L\in\mathcal L_m(G_{\ell,t})\), and write
\(\mathcal U(L)=\bigcup_v L(v)\).
For \(c\in\mathcal U(L)\), let \(d_c\) be the number of
lists containing \(c\). Admissibility implies that the
vertices whose lists contain \(c\) form a clique, so
\(d_c\le2t\). Counting vertex-color incidences gives
\(\ell tm =\sum_{c\in\mathcal U(L)}d_c \le2t\,|\mathcal U(L)|\).
Consequently,
\(|\mathcal U(L)| \ge\left\lceil\frac{\ell m}{2}\right\rceil\).

We construct an admissible assignment attaining this bound.
Choose pairwise disjoint palettes
\(P_0,\ldots,P_{\ell-1}\), with their indices also taken
modulo \(\ell\).
If \(m\) is even, let \(|P_i|=m/2\) for every \(i\), and set
\(L(v)=P_{i-1}\cup P_i \qquad (v\in V_i)\).
Every list has size \(m\), and the union of the lists has
size \(\ell m/2\).

Suppose that \(m\) is odd, and put \(a=(m-1)/2\).
Choose the palette sizes so that
\[
 |P_i|=
 \begin{cases}
 a+1, & i\in\{0,2,\ldots,\ell-3\},\\
 a,   & \text{otherwise}.
 \end{cases}
\]
Let \(\gamma\) be a single color outside
\(\bigcup_{i=0}^{\ell-1}P_i\), and define
\[
 L(v)=
 \begin{cases}
 P_{i-1}\cup P_i,
     & v\in V_i,\quad 0\le i\le\ell-2,\\
 P_{\ell-2}\cup P_{\ell-1}\cup\{\gamma\},
     & v\in V_{\ell-1}.
 \end{cases}
\]
For \(0\le i\le\ell-2\), exactly one of
\(P_{i-1}\) and \(P_i\) has size \(a+1\), so these lists
have size \(2a+1=m\).
Both \(P_{\ell-2}\) and \(P_{\ell-1}\) have size \(a\);
adding the same color \(\gamma\) at every vertex of
\(V_{\ell-1}\) gives lists of size \(m\) there as well.
There are \(k\) palettes of size \(a+1\) and \(k+1\)
palettes of size \(a\). Hence
\[
 |\mathcal U(L)|
 =k(a+1)+(k+1)a+1
 =\ell a+k+1
 =\left\lceil\frac{\ell m}{2}\right\rceil.
\]

In either construction, each color in \(P_i\) occurs
precisely on the clique \(V_i\cup V_{i+1}\).
In the odd case, the additional color \(\gamma\) occurs
precisely on the clique \(V_{\ell-1}\).
Thus the support of every color is a clique, and the
assignment is admissible.
This proves the asserted minimum for the union of the lists.
Every proper coloring uses exactly \(\ell t\) distinct
colors, so
\(\mathcal C_{G_{\ell,t}}(L) \subseteq \{S\subseteq\mathcal U(L):|S|=\ell t\}\).
It follows that
\[
 f'(G_{\ell,t},m)
 \le |\mathcal C_{G_{\ell,t}}(L)|
 \le \binom{\lceil\ell m/2\rceil}{\ell t}.
\]

It remains to determine the clique-partition value.
The part sizes of any clique partition of \(G_{\ell,t}\)
are positive integers with sum \(\ell t\), each at most
\(2t\). We first minimize
\(\prod_j\binom{m}{a_j}\) over all integer partitions
with these restrictions. Since every clique partition
provides such an integer partition, this minimum is a
lower bound for \(w(G_{\ell,t},m)\).

Consider an integer partition minimizing this product.
If it contains two parts \(u,v\) with
\(1\le v\le u<2t\), replace them by \(u+1,v-1\), omitting
the zero part when \(v=1\).
The sum of the parts is unchanged, and every remaining
part still lies between \(1\) and \(2t\).
With \(\binom{m}{0}=1\), the ratio of the new product
to the old product is
\[
 \frac{\binom{m}{u+1}\binom{m}{v-1}}
      {\binom{m}{u}\binom{m}{v}}
 =
 \frac{(m-u)v}{(u+1)(m-v+1)}
 <1.
\]
Indeed, the denominator exceeds the numerator by
\((u+1)(m-v+1)-(m-u)v =(m+1)(u+1-v)>0\).
This contradicts minimality. Hence at most one part
is smaller than \(2t\).
Since
\(\ell t=k(2t)+t, \qquad 0<t<2t\),
the minimizing integer partition consists of exactly
\(k\) parts of size \(2t\) and one part of size \(t\).
Therefore
\(w(G_{\ell,t},m) \ge \binom{m}{2t}^{k}\binom{m}{t}\).

Finally,
\[
 \{V_0\cup V_1,\ V_2\cup V_3,\ \ldots,\
   V_{\ell-3}\cup V_{\ell-2},\ V_{\ell-1}\}
\]
is a clique partition of \(G_{\ell,t}\) with precisely
these part sizes. It attains the preceding lower bound,
and hence
\(w(G_{\ell,t},m) =\binom{m}{2t}^{k}\binom{m}{t}\).
\end{proof}

For the \(20\)-vertex graph, the two quantities in the
lemma are explicit integers. For the fixed \(25\)-vertex
graph, their ratio depends on the list size; a factor-by-factor
estimate bounds it uniformly away from one and gives
the strict comparison throughout the required range.

\begin{proof}[Proof of Theorem~\ref{thm:counterexamples}]
Applying Lemma~\ref{lem:palettes} with
\(\ell=5\), \(t=4\), and \(m=20\), we obtain
\[
 \begin{aligned}
 f'(G_{5,4},20)
 &\le \binom{50}{20}
  =47\,129\,212\,243\,960\\
 &<76\,882\,596\,160\,500
  =\binom{20}{8}^{2}\binom{20}{4}
  =w(G_{5,4}).
 \end{aligned}
\]
Since \(G_{5,4}\) has \(20\) vertices,
Lemma~\ref{lem:join} gives
\(f(G_{5,4})=f'(G_{5,4},20)\), proving the first assertion.

For the second assertion, fix an integer \(m\ge25\), and put
\(q=\lceil5m/2\rceil\).
Lemma~\ref{lem:palettes}, applied with \(\ell=t=5\), gives
\(f'(G_{5,5},m)\le\binom{q}{25}, \qquad w(G_{5,5},m)=\binom{m}{10}^{2}\binom{m}{5}\).
It therefore suffices to prove
\(\binom{m}{10}^{2}\binom{m}{5}>\binom{q}{25}\).
Let \(a_0\le\cdots\le a_{24}\) be the sequence obtained
by sorting two copies of \(0,\ldots,9\) and one copy
of \(0,\ldots,4\). We claim that
\(5a_j\le2j-1 \qquad (1\le j\le24)\).
If \(a_j=0\), the inequality is immediate.
For \(1\le u\le4\), each smaller value occurs three times,
so the first occurrence of \(u\) has index \(3u\).
Thus, whenever \(a_j=u\),
\(2j-1\ge6u-1\ge5u=5a_j\).
For \(5\le u\le9\), there are \(15\) entries smaller than
\(5\), and each value from \(5\) to \(u-1\) occurs twice.
The first occurrence of \(u\) therefore has index
\(15+2(u-5)=2u+5\). Hence
\(2j-1\ge4u+9\ge5u=5a_j\),
where the last inequality uses \(u\le9\).
This proves the claim.

Since \(m\) is an integer,
\(q\le(5m+1)/2\). Consequently, for \(1\le j\le24\),
\(q-j \le\frac{5m+1-2j}{2} \le\frac52(m-a_j)\).
All these quantities are positive, because \(q\ge63\)
and \(a_j\le9<m\). Therefore
\(\frac{m-a_j}{q-j}\ge\frac25 \qquad (1\le j\le24)\).
For the remaining factor, \(a_0=0\) gives
\[
 \frac{m-a_0}{q}
 =\frac{m}{q}
 \ge\frac{2m}{5m+1}
 =\frac25\,\frac{5m}{5m+1}.
\]

Expanding the binomial coefficients into falling factorials
and using the multiplicities defining the sequence
\((a_j)_{j=0}^{24}\), we obtain
\[
 \begin{aligned}
 \frac{\binom{m}{10}^{2}\binom{m}{5}}
      {\binom{q}{25}}
 &=
 \frac{25!}{(10!)^2\,5!}
 \prod_{j=0}^{24}\frac{m-a_j}{q-j}\\
 &\ge
 \frac{25!}{(10!)^2\,5!}
 \left(\frac25\right)^{25}
 \frac{5m}{5m+1}.
 \end{aligned}
\]
Put
\(C=\frac{25!}{(10!)^2\,5!} \left(\frac25\right)^{25}\).
An exact calculation gives
\[
 C=
 \frac{65\,874\,639\,917\,678\,592}
      {59\,604\,644\,775\,390\,625}
 >\frac{11}{10}.
\]
Moreover, \(m\ge25\) implies
\(5m/(5m+1)\ge125/126\). Hence
\[
 \frac{\binom{m}{10}^{2}\binom{m}{5}}
      {\binom{q}{25}}
 \ge C\,\frac{5m}{5m+1}
 >\frac{11}{10}\cdot\frac{125}{126}
 =\frac{275}{252}>1.
\]
It follows that
\[
 f'(G_{5,5},m)
 \le\binom{q}{25}
 <\binom{m}{10}^{2}\binom{m}{5}
 =w(G_{5,5},m).
\]
Since \(m\ge25\) was arbitrary, the second assertion follows.
\end{proof}

The binomial upper bound also counts subsets that cannot
be realized by a coloring. Appendix~\ref{sec:finite-count}
identifies and counts the realizable subsets of the
\(50\)-color assignment, giving a sharper upper bound
for \(f(G_{5,4})\).

\section{The triangle-free case}\label{sec:trianglefree}

The proof of Theorem~\ref{thm:trianglefree} will also give
the following quantitative estimate. Its correction term
provides the strictness needed to determine which list
assignments attain the minimum.

\begin{proposition}\label{prop:triangle-quant}
Let \(G\) be a triangle-free graph, and put
\(n=|V(G)|\ge2\) and \(\nu=\nu(G)\).
Let \(m\ge n\) be an integer, and let
\(L\in\mathcal L_m(G)\).
For \(e=uv\in E(G)\), set
\(p_e=p_{uv}=\frac{|L(u)\cap L(v)|}{m}, \qquad R=\nu-\sum_{e\in E(G)}p_e^2\).
Then \(R\ge0\), and
\[
 N_L\ge
 \binom{m}{2}^{\nu}m^{n-2\nu}
 +\frac{m^n}{40\cdot2^\nu}\,R.
\]
\end{proposition}

\subsection{Color sets and orientations}

Fix \(G,m,L\) as in the proposition, and write
\(V=V(G)\), \(E=E(G)\), and \(\nu=\nu(G)\).
Let \(\mathcal U=\bigcup_{v\in V}L(v)\).
The support of every color is a clique, so triangle-freeness
implies that each color occurs in either one or two lists.
We call these colors private and shared, respectively.
In particular, the shared colors in a fixed list \(L(v)\)
are partitioned according to their other endpoint. Hence
\[
 \sum_{\substack{e\in E\\e\ni v}}p_e\le1
 \qquad (v\in V),
\]
and \(p=(p_e)_{e\in E}\) belongs to \(\mathcal P_f(G)\).

We first verify that \(R\ge0\).
By Lemma~\ref{lem:fractional}, the positive support of an
extreme point \(x\in\mathcal P_f(G)\) consists of disjoint
edges of weight one and vertex-disjoint odd cycles of
weight \(1/2\).
If there are \(a\) such edges and the cycle lengths are
\(\ell_1,\ldots,\ell_k\), then
\[
 \sum_{e\in E}x_e^2
 =a+\sum_{i=1}^k\frac{\ell_i}{4}
 \le a+\sum_{i=1}^k\frac{\ell_i-1}{2}
 \le\nu.
\]
The last inequality follows by taking the weight-one edges
and a maximum matching from each of the disjoint cycles.
Since every fractional matching is a convex combination
of extreme points, convexity of the squared Euclidean
norm gives
\(\sum_{e\in E}p_e^2\le\nu\).
Thus \(R\ge0\).

We now associate a multigraph \(H=H_L\) with the assignment.
Its vertex set is \(V\), and each color \(c\in\mathcal U\)
corresponds to one edge \(e_c\).
A private color at \(v\) gives a loop at \(v\), while a
color shared by \(u\) and \(v\) gives an edge joining them.
Distinct colors give distinct labeled edges, even when
their endpoints coincide.
For \(S\subseteq\mathcal U\), let
\(H_S=(V,\{e_c:c\in S\})\).
Loops are regarded as cycles of length one, and two
distinct parallel edges form a cycle of length two.
In applying Lemma~\ref{lem:cycle-orientations}, the
representatives are assigned to vertices and the labeled
edges are colors. The direction of a selected edge records
which incident vertex receives its color. With this
convention, both the criterion for a realizable image and
the number of maps realizing it can be read from its
connected components.

A set \(S\) belongs to \(\mathcal C_G(L)\) if and only if
every connected component of \(H_S\) is unicyclic.
Indeed, a proper \(L\)-coloring with image \(S\) assigns
each color in \(S\) to exactly one of its incident vertices.
Orienting its edge away from that vertex gives an
orientation of \(H_S\) with outdegree one at every vertex.
Each connected component consequently has as many edges
as vertices and is therefore unicyclic.
Conversely, orienting the unique cycle of each component
cyclically and directing every remaining edge towards
that cycle gives outdegree one at every vertex.
Assigning to each vertex the color of its outgoing edge
then yields a proper \(L\)-coloring with image \(S\).

For \(S\in\mathcal C_G(L)\), let \(\kappa(S)\) be the
number of components of \(H_S\) whose cycle has length
at least two.
The directions of all edges outside the cycles are forced.
A loop has one orientation, whereas every nonloop cycle
has two. This includes a cycle of two parallel edges,
whose distinct color labels distinguish the two orientations.
Thus
\[
 \bigl|\{\varphi:
       \varphi\text{ is a proper }L\text{-coloring},
       \ \varphi(V)=S\}\bigr|
 =2^{\kappa(S)}.
\]
This counts the different coloring maps in the multigraph
extension of Lemma~\ref{lem:cycle-orientations};
their common image contributes only once to the color-set
count. Selecting the same shared color at both endpoints
would instead use one labeled edge twice. It must therefore
be distinguished from selecting two different parallel
color edges when we define the collision event below.

Choose a color \(\xi_v\) independently and uniformly from
\(L(v)\) at each vertex \(v\).
The sample space
\(\Omega=\prod_{v\in V}L(v)\) consists of \(m^n\)
equally likely labeled color choices.
Each outcome induces a map \(F:V\to V\):
a private color chosen at \(v\) gives \(F(v)=v\), while
a color shared with \(u\ne v\) gives \(F(v)=u\).
For every edge \(uv\in E\),
\(\PP(F(u)=v)=p_{uv}=p_{vu}=\PP(F(v)=u)\).
Let \(B\) be the event that all chosen colors are distinct.
By admissibility, \(B\) is precisely the event that
\(\xi=(\xi_v)_{v\in V}\) is a proper \(L\)-coloring.
Let \(C\) count the directed cycles of \(F\) of length
at least two.
Distinct directed cycles are vertex-disjoint, and selecting
one underlying edge from each such cycle gives a matching
in \(G\). Consequently,
\(C\le\nu\)
for every outcome, including outcomes outside \(B\).

On \(B\), an outcome with image \(S\) satisfies
\(C=\kappa(S)\).
Grouping the outcomes in \(B\) according to their images
therefore gives
\[
 \begin{aligned}
 \E\!\left[2^{-C}\mathbf{1}_B\right]
 &=
 \frac{1}{m^n}
 \sum_{S\in\mathcal C_G(L)}
       2^{-\kappa(S)}\,2^{\kappa(S)}\\
 &=\frac{N_L}{m^n}.
 \end{aligned}
\]
For subsequent estimates, put
\(s=\sum_{e\in E}p_e, \qquad u_2=\sum_{e\in E}p_e^2, \qquad D=\nu-\E C\).
The fractional matching constraints imply \(s\le n/2\),
and the cycle bound gives \(D\ge0\).

Since \(\nu-C\) is a nonnegative integer,
\(2^{\nu-C}\ge1+\nu-C\). Hence
\begin{equation}\label{eq:fibers}
 \frac{2^\nu N_L}{m^n}
 =
 \E\!\left[2^{\nu-C}\mathbf{1}_B\right]
 \ge
 \PP(B)+\E\!\left[(\nu-C)\mathbf{1}_B\right].
\end{equation}
Using
\(\binom{m}{2}^{\nu}m^{n-2\nu} =\frac{m^n}{2^\nu}\left(1-\frac1m\right)^\nu\),
we see that it suffices to establish
\[
 \PP(B)+\E\!\left[(\nu-C)\mathbf{1}_B\right]
 \ge
 \left(1-\frac1m\right)^\nu+\frac{R}{40}.
\]
The following estimates control the collision event and
the cycle deficit jointly.

\subsection{Directed cycles and collisions}

The first estimate bounds the probability that a prescribed
arc lies on a directed cycle of length at least three.
Following the orbit from its head leads to an escape event
for a reversible Markov chain. After deleting the underlying
edge, Lemma~\ref{lem:dirichlet} bounds that escape probability
by the energy of a test function. The lemma is stated for
an arbitrary simple graph; triangle-freeness will enter
the subsequent estimate for the total number of cycles.

\begin{lemma}\label{lem:conductance}
Let \(G=(V,E)\) be a finite simple graph.
Assign weights \(p_{xy}=p_{yx}\ge0\) to its edges, with
\(\sum_{y:xy\in E}p_{xy}\le1 \qquad (x\in V)\).
Let \(F:V\to V\) be a random map whose values are chosen
independently, with
\[
 \PP(F(x)=y)=p_{xy}\quad(xy\in E),
 \qquad
 \PP(F(x)=x)=1-\sum_{y:xy\in E}p_{xy}.
\]
For an edge \(e=uv\), put \(p_e=p_{uv}\), and let
\(A_{v\to u}\) be the event that \(F(v)=u\) and the arc
\(v\to u\) belongs to a directed cycle of length at least
three. Then
\[
 \PP(A_{v\to u})\le\frac{p_e(1-p_e)}{2}.
\]
\end{lemma}

\begin{proof}
Put \(a=p_e\).
If \(a=0\), the event \(A_{v\to u}\) is impossible.
If \(a=1\), symmetry and the degree constraints imply
\(F(u)=v\) and \(F(v)=u\) almost surely, so these vertices
form a directed cycle of length two.
The assertion therefore holds in both cases.
Henceforth assume \(0<a<1\).

Delete \(e\), and define symmetric conductances on \(V\) by
\[
 c_{xy}=
 \begin{cases}
 p_{xy}, & xy\in E\setminus\{e\},\\
 0, & \text{otherwise}.
 \end{cases}
\]
In particular, \(c_{xx}=0\).
Let \(P'\) be the transition matrix given by
\[
 P'_{xy}=c_{xy}\quad(x\ne y),
 \qquad
 P'_{xx}=1-\sum_{y\in V}c_{xy}.
\]
Thus the deleted transition mass is added to the holding
probabilities at \(u\) and \(v\).
For the Markov chain \(X=(X_j)_{j\ge0}\) with transition
matrix \(P'\), write \(\PP'_x\) for its law started at \(x\),
and define
\[
 \tau_x=\inf\{j\ge0:X_j=x\},
 \qquad
 \tau_x^+=\inf\{j\ge1:X_j=x\}.
\]
As usual, the infimum of the empty set is infinity.
Put
\(h(x)=\PP'_x(\tau_v<\tau_u)\).
Then \(h(u)=0\), \(h(v)=1\), and \(h\) is harmonic
on \(V\setminus\{u,v\}\).
On a component of the positive-conductance graph meeting
neither boundary vertex, \(h\) is identically zero.

Conditional on \(F(v)=u\), the event \(A_{v\to u}\) occurs
precisely when the forward orbit from \(u\) follows a simple
\(u\)-\(v\) path in \(G-e\).
Let \(\Pi_{uv}\) be the finite family of all such paths.
For a path \(\pi=(u=x_0,x_1,\ldots,x_k=v)\), the required
arrows are chosen at distinct vertices different from \(v\).
Their conditional probability is therefore
\(\prod_{i=0}^{k-1}p_{x_i x_{i+1}}\), which is also the
probability that the Markov chain with transition matrix
\(P'\) follows this path. Denote this common path weight
by \(w(\pi)\).

To make the comparison of the two probability spaces
explicit, let \(E'_\pi=\{X_i=x_i:0\le i\le k\}\)
be the trajectory event under \(\PP'_u\) associated
with \(\pi\).
In the random-map space, the events prescribing the
arrows along the distinct paths in \(\Pi_{uv}\) are
disjoint under the conditional law given \(F(v)=u\),
and their union is the conditional cycle event.
Consequently,
\[
 \PP(A_{v\to u}\mid F(v)=u)
 =\sum_{\pi\in\Pi_{uv}}w(\pi)
 =\sum_{\pi\in\Pi_{uv}}\PP'_u(E'_\pi).
\]

In the Markov-chain space, the events corresponding to
distinct simple \(u\)-\(v\) paths are also disjoint:
each specified trajectory reaches \(v\) for the first
time at its final step. Moreover, each such event is
contained in \(\{\tau_v<\tau_u^+\}\) for the Markov
chain started at \(u\), because a simple path does not
revisit its initial vertex before reaching its endpoint.
The chain can additionally reach its endpoint after
holding steps or repeated visits to intermediate vertices,
which are not part of these simple-path events.
Thus comparison is made through their probabilities,
and we obtain
\[
 \PP(A_{v\to u}\mid F(v)=u)
 \le \PP'_u(\tau_v<\tau_u^+)
 =\sum_{x\in V}c_{ux}h(x).
\]
Denote the last sum by \(J\).

For a function \(g:V\to\mathbb R\), define
\[
 \mathcal E(g)=
 \sum_{xy\in E\setminus\{e\}}
       c_{xy}\bigl(g(x)-g(y)\bigr)^2,
\]
where each undirected edge is counted once.
Lemma~\ref{lem:dirichlet} gives
\[
 \mathcal E(h)
 =
 \min\{\mathcal E(g):g(u)=0,\ g(v)=1\}.
\]
To identify this energy with \(J\), put
\(b(x)=\sum_{y\in V}c_{xy}\bigl(h(x)-h(y)\bigr)\).
Harmonicity gives \(b(x)=0\) for \(x\notin\{u,v\}\).
Moreover, \(b(u)=-J\), and symmetry implies
\(\sum_{x\in V}b(x)=0\). Hence \(b(v)=J\).
Expanding the energy now yields
\(\mathcal E(h) =\sum_{x\in V}h(x)b(x) =b(v)=J\).

Consider the test function \(g_0\) with
\(g_0(u)=0\), \(g_0(v)=1\), and \(g_0(x)=1/2\)
for \(x\notin\{u,v\}\).
Since \(c_{uv}=0\), only edges incident with \(u\) or \(v\)
contribute to its energy. Thus
\[
 \begin{aligned}
 J
 =\mathcal E(h)
 &\le\mathcal E(g_0)\\
 &=\frac14
   \left(\sum_{x\in V}c_{ux}
        +\sum_{x\in V}c_{vx}\right)\\
 &\le\frac14\bigl((1-a)+(1-a)\bigr)
 =\frac{1-a}{2}.
 \end{aligned}
\]
Finally, \(\PP(F(v)=u)=a\), so
\[
 \PP(A_{v\to u})
 =a\,\PP(A_{v\to u}\mid F(v)=u)
 \le aJ
 \le\frac{a(1-a)}2.
\]
This proves the lemma.
\end{proof}

\begin{lemma}\label{lem:cycles-collisions}
Let \(G\) be a triangle-free graph on \(n\ge2\) vertices,
let \(m\ge n\), and let \(L\in\mathcal L_m(G)\).
For the associated random map, retain the notation
\(B,C,\nu,s,u_2,D\) introduced above. Then
\[
 \E C\le u_2+\frac{s-u_2}{4},
 \qquad
 \E[(\nu-C)\mathbf 1_B]
 \ge \left(1-\frac{s}{m}\right)D
 \ge \frac D2.
\]
\end{lemma}

\begin{proof}
All probabilities and expectations refer to the independent
uniform choices from the lists. All edge sums are over \(E(G)\).
Let \(C_2\) denote the number of directed two-cycles.
For each edge \(e=uv\), the probability that the random map
contains the two-cycle \(u\to v\to u\) is \(p_e^2\).
Hence
\(\E C_2=\sum_e p_e^2=u_2\).
For an oriented edge \(x\to y\), let \(A_{x\to y}\) be the
event that this arc belongs to a directed cycle of length
at least three. Since \(G\) is triangle-free, every such
cycle has length at least four. Lemma~\ref{lem:conductance}
therefore gives
\[
 \begin{aligned}
 4\E(C-C_2)
 &\le
 \sum_{uv\in E(G)}
 \bigl(\PP(A_{u\to v})+\PP(A_{v\to u})\bigr)\\
 &\le \sum_e p_e(1-p_e)
 =s-u_2.
 \end{aligned}
\]
Consequently,
\begin{equation}\label{eq:defect}
 \E C\le u_2+\frac{s-u_2}{4}.
\end{equation}

For an edge \(e=uv\), let \(T_e\) be the number of directed
cycles of length at least two that meet \(\{u,v\}\).
A cycle containing both endpoints is counted once.
Since the directed cycles of a map are vertex-disjoint,
\(T_e\in\{0,1,2\}\): at most one cycle can pass through
each endpoint. We next prove
\begin{equation}\label{eq:touch}
 \sum_e p_e\E T_e\le s.
\end{equation}
Write
\[
 d_p(x)=\sum_{e\ni x}p_e\le1,
 \qquad
 \eta(W)=
 \sum_{\substack{e\in E(G)\\e\cap W\ne\varnothing}}p_e
 \quad (W\subseteq V(G)).
\]
Thus
\(\eta(W) =\sum_{x\in W}d_p(x) -\sum_{e\in E(G[W])}p_e\).
For a realization of the random map, let \(\mathcal K\)
be its family of directed cycles of length at least two.
Double counting gives
\[
 \sum_e p_eT_e
 =\sum_{K\in\mathcal K}\eta(V(K)).
\]

If \(K\) is a two-cycle on \(e=uv\), then
\(\eta(V(K)) =d_p(u)+d_p(v)-p_e \le2-p_e\).
The expected contribution of all two-cycles is therefore
at most
\(\sum_e p_e^2(2-p_e)\).
If \(K\) is a longer cycle, let \(E(K)\) denote its set
of underlying undirected edges. Since
\(E(K)\subseteq E(G[V(K)])\), we have
\[
 \eta(V(K))
 \le |V(K)|-\sum_{e\in E(K)}p_e
 =\sum_{e\in E(K)}(1-p_e).
\]
By Lemma~\ref{lem:conductance}, the expected contribution
of all longer cycles is at most
\[
 \begin{aligned}
 &\sum_{uv\in E(G)}(1-p_{uv})
   \bigl(\PP(A_{u\to v})+\PP(A_{v\to u})\bigr)\\
 &\hspace{3em}\le
 \sum_e p_e(1-p_e)^2.
 \end{aligned}
\]
Combining these estimates yields
\[
 \sum_e p_e\E T_e
 \le
 \sum_e\bigl(p_e^2(2-p_e)+p_e(1-p_e)^2\bigr)
 =\sum_e p_e=s,
\]
which proves \eqref{eq:touch}.

Let \(\xi_v\) denote the color chosen at \(v\).
For \(e=uv\), define the collision event
\(F_e=\{\xi_u=\xi_v\}\).
Independence and uniformity give
\(\PP(F_e) =\frac{|L(u)\cap L(v)|}{m^2} =\frac{p_e}{m}\).
Let \(C_e^\circ\) be the number of directed cycles of
length at least two that avoid both \(u\) and \(v\).
Then
\(C=T_e+C_e^\circ\).
The variable \(C_e^\circ\) depends only on the choices
at vertices outside \(\{u,v\}\), so it is independent
of \(F_e\). Indeed, the collision event depends only
on the choices at the two endpoints, and the vertex
choices are independent. Conditioning on that event
therefore leaves the joint law of all choices outside
the endpoints unchanged, and in particular leaves
the expectation of \(C_e^\circ\) unchanged.
On \(F_e\), the endpoints form the two-cycle
\(u\to v\to u\), and hence \(T_e=1\).
Thus, whenever \(p_e>0\),
\[
 \begin{aligned}
 \E[\nu-C\mid F_e]
 &=\nu-1-\E C_e^\circ\\
 &=\nu-\E C+\E T_e-1\\
 &=D+\E T_e-1.
 \end{aligned}
\]

Every repeated color occurs at the endpoints of an edge,
so
\(B^c=\bigcup_e F_e\).
Moreover, \(C\le\nu\) for every realization: choosing one
underlying edge from each nontrivial directed cycle gives
a matching. Put \(Z=\nu-C\ge0\), so that \(\E Z=D\).
The pointwise inequality
\(Z\mathbf 1_{B^c}\le\sum_e Z\mathbf 1_{F_e}\)
uses only the nonnegativity of the deficit, so it remains
valid even when the collision events on different edges
are dependent. Combining this weighted union bound
with the conditional identity just proved and
\eqref{eq:touch}, we obtain
\begin{equation}\label{eq:noncollision}
 \begin{aligned}
 \E[Z\mathbf 1_B]
 &=D-\E[Z\mathbf 1_{B^c}]\\
 &\ge
 D-\sum_{\substack{e\in E(G)\\p_e>0}}
       \frac{p_e}{m}\bigl(D+\E T_e-1\bigr)\\
 &=\left(1-\frac{s}{m}\right)D
   -\frac1m\left(\sum_e p_e\E T_e-s\right)\\
 &\ge\left(1-\frac{s}{m}\right)D.
 \end{aligned}
\end{equation}
Finally,
\(2s=\sum_{v\in V(G)}d_p(v)\le n\le m\).
Since \(D\ge0\), the last expression in
\eqref{eq:noncollision} is at least \(D/2\).
\end{proof}

\subsection{Matching bounds and equality cases}

The fractional matching decomposition has two uses below.
It bounds the cycle deficit in terms of the squared overlaps
and provides a matching on which collision events are
independent. These are the inputs for comparing the
color-set count with the matching construction.

\begin{proof}[Proof of Proposition~\ref{prop:triangle-quant}]
Assume first that \(n\ge2\). Fix \(L\in\mathcal L_m(G)\),
write \(N_L=|\mathcal C_G(L)|\), and retain the notation
\(\nu,p,s,u_2,D,R\) introduced above.
Since every color occurs in at most two lists, the sets
\(L(v)\cap L(u)\), as \(u\) ranges over the neighbors of \(v\),
are pairwise disjoint. Consequently,
\(\sum_{e\ni v}p_e\le1 \qquad\text{for every }v\in V(G)\),
so \(p\in\mathcal P_f(G)\).
Express \(p\) as a convex combination of extreme points
of \(\mathcal P_f(G)\), and let \(X\) be the corresponding
random extreme point. Thus \(\E X=p\).

By Lemma~\ref{lem:fractional}, the positive entries of \(X\)
are \(1\) on a matching and \(1/2\) on vertex-disjoint odd
cycles, all disjoint from the matching edges.
Let \(k\) be the number of \(1\)-edges, and let
\(K_1,\ldots,K_b\) be these odd cycles, with lengths
\(\ell_1,\ldots,\ell_b\). Since \(G\) is triangle-free,
each \(\ell_i\ge5\). Taking a maximum matching in each
\(K_i\), together with the \(1\)-edges, gives
\(\nu\ge k+\sum_{i=1}^b\frac{\ell_i-1}{2}\).
It follows that
\[
 \begin{aligned}
 \nu-\sum_e X_e^2
 &\ge
 \sum_{i=1}^b\frac{\ell_i-2}{4}
 \ge\frac34 b,\\
 \sum_e\left(X_e+\frac23X_e^2\right)
 &=\frac53k+\frac23\sum_{i=1}^b\ell_i\\
 &\le
 \frac53\left(k+\sum_{i=1}^b\frac{\ell_i-1}{2}\right)
 \le\frac53\nu.
 \end{aligned}
\]
Here the second estimate uses
\(2\ell_i/3\le5(\ell_i-1)/6\), which holds because
\(\ell_i\ge5\).
Taking expectations and applying Jensen's inequality to
the squared Euclidean norm, we obtain
\[
 u_2=\|p\|_2^2
 \le\E\|X\|_2^2
 \le\nu-\frac34\E b,
 \qquad
 s+\frac23u_2\le\frac53\nu.
\]
Therefore,
\(R\ge0,\qquad \E b\le\frac43R,\qquad s-\nu\le\frac23R\).
In particular,
\(s-u_2=(s-\nu)+R\le\frac53R\).
By \eqref{eq:defect},
\[
 D=\nu-\E C
 \ge R-\frac{s-u_2}{4}
 \ge\frac7{12}R.
\]

We next find a matching of large \(p\)-weight.
Conditional on \(X\), retain every \(1\)-edge and choose,
independently in each \(K_i\), a uniformly random maximum
matching. These choices together form a matching
\(\widehat M\).
By symmetry, each edge of \(K_i\) is selected with
probability \((\ell_i-1)/(2\ell_i)\). Hence
\[
 \E\left[\sum_{e\in\widehat M}p_e\,\middle|\,X\right]
 =
 p\cdot X
 -\sum_{i=1}^b\frac1{2\ell_i}
                  \sum_{e\in E(K_i)}p_e.
\]
Summing the fractional matching constraints over
\(V(K_i)\) gives
\[
 2\sum_{e\in E(K_i)}p_e
 \le\sum_{v\in V(K_i)}\sum_{e\ni v}p_e
 \le\ell_i.
\]
Thus each cycle contributes a loss of at most \(1/4\).
Taking expectations yields
\[
 \E\sum_{e\in\widehat M}p_e
 \ge p\cdot\E X-\frac14\E b
 =u_2-\frac14\E b
 \ge\nu-\frac43R.
\]
There is therefore a matching \(M\) whose weight
\(W=\sum_{e\in M}p_e\)
satisfies \(W\ge\nu-4R/3\). Also \(W\le|M|\le\nu\),
since \(0\le p_e\le1\).

Put \(a=1-1/m\). The collision events on the edges of
\(M\) depend on disjoint sets of vertex choices and are
therefore independent. A union bound for the collision
events outside \(M\) gives
\[
 \PP(B)
 \ge\prod_{e\in M}(1-p_e/m)-\frac{s-W}{m}.
\]
By concavity,
\(\log(1-x/m)\ge x\log(1-1/m) \qquad (0\le x\le1)\).
Consequently,
\(\prod_{e\in M}(1-p_e/m)\ge a^W\).
Moreover, Bernoulli's inequality and
\(\nu\le n/2\le m/2\) imply
\(a^\nu\ge1-\frac{\nu}{m}\ge\frac12\).
Since \(W\le\nu\) and \(-\log a\ge1/m\), we have
\[
 \begin{aligned}
 a^W-a^\nu
 &=a^\nu
   \left(\exp\bigl((\nu-W)(-\log a)\bigr)-1\right)\\
 &\ge a^\nu(\nu-W)(-\log a)
 \ge\frac{\nu-W}{2m}.
 \end{aligned}
\]
It follows that
\begin{equation}\label{eq:coll}
 \begin{aligned}
 \PP(B)
 &\ge
 a^\nu-\frac{s-\nu+(\nu-W)/2}{m}\\
 &\ge a^\nu-\frac{4R}{3m},
 \end{aligned}
\end{equation}
where the last step uses
\(s-\nu\le2R/3\) and \(\nu-W\le4R/3\).

Suppose first that \(m\ge5\). Combining
\eqref{eq:fibers}, \eqref{eq:noncollision},
and \eqref{eq:coll}, we obtain
\[
 \begin{aligned}
 \frac{2^\nu N_L}{m^n}
 &\ge \PP(B)+\E[(\nu-C)\mathbf 1_B]\\
 &\ge a^\nu-\frac{4R}{3m}+\frac D2\\
 &\ge a^\nu+
       \left(\frac7{24}-\frac4{3m}\right)R\\
 &\ge a^\nu+\frac{R}{40}.
 \end{aligned}
\]
The last inequality follows from
\(\frac7{24}-\frac4{3m} \ge\frac7{24}-\frac4{15} =\frac1{40}\).

It remains to consider \(2\le m\le4\).
Then \(n\le4\), so the triangle-free graph \(G\) is
bipartite. By Lemma~\ref{lem:fractional}, every extreme
point \(X\) is the indicator of a matching. Hence
\(s=\E|X|\le\nu\).
Also,
\(\E\sum_e p_eX_e=u_2\),
so we may choose a matching \(M\) with weight
\(W\ge u_2=\nu-R\).
Since \(s-u_2\le R\), \eqref{eq:defect} gives
\(D\ge3R/4\).
The first inequality in \eqref{eq:coll}, applied to this
matching, yields
\[
 \PP(B)
 \ge a^\nu-\frac{s-\nu+(\nu-W)/2}{m}
 \ge a^\nu-\frac{R}{2m}.
\]
Using \eqref{eq:fibers} and \eqref{eq:noncollision} again,
we obtain
\[
 \begin{aligned}
 \frac{2^\nu N_L}{m^n}
 &\ge a^\nu-\frac{R}{2m}+\frac D2\\
 &\ge a^\nu+
       \left(\frac38-\frac1{2m}\right)R\\
 &\ge a^\nu+\frac R8.
 \end{aligned}
\]
Thus, for every \(m\ge n\ge2\),
\[
 \frac{2^\nu N_L}{m^n}
 \ge a^\nu+\frac{R}{40}.
\]
Since
\(\frac{m^n}{2^\nu}a^\nu =\binom m2^\nu m^{n-2\nu}\),
this proves Proposition~\ref{prop:triangle-quant}.
\end{proof}

The quantitative lower bound and the clique-partition
upper bound determine the extremal value. For the equality
statement, the remaining task is to identify assignments
whose correction term vanishes. Strict convexity in the
extreme-point decomposition then determines the overlap
vector and the supporting maximum matching.

\begin{proof}[Proof of Theorem~\ref{thm:trianglefree}]
Assume first that \(n=|V(G)|\ge2\), and put \(\nu=\nu(G)\).
Proposition~\ref{prop:triangle-quant} applies to every
\(L\in\mathcal L_m(G)\).
Dropping the nonnegative correction term and minimizing
over \(L\in\mathcal L_m(G)\), we obtain
\[
 f'(G,m)\ge
 \binom m2^\nu m^{n-2\nu}=w(G,m).
\]
The reverse inequality follows from the clique-partition
construction, so equality holds.

We finish by determining all minimizing assignments.
If \(N_L=w(G,m)\), the quantitative bound forces \(R=0\).
For this assignment, let
\(p_e=|L(u)\cap L(v)|/m\) for \(e=uv\), and choose a random
extreme point \(X\) of \(\mathcal P_f(G)\) with \(\E X=p\).
Let \(b\) count its odd cycles of weight \(1/2\).
The extreme-point estimates in the preceding proof give
\(\E\|X\|_2^2\le\nu\) and
\(\nu-\|X\|_2^2\ge3b/4\) for every realization.
Since \(\E X=p\), we have
\[
 \begin{aligned}
 0
 &\le\E\|X-p\|_2^2\\
 &=\E\|X\|_2^2-\|p\|_2^2\\
 &\le\nu-u_2
 =0.
 \end{aligned}
\]
Thus \(X=p\) almost surely.
The inequality
\(\nu-\|X\|_2^2\ge3b/4\) then implies \(b=0\).
Consequently, \(p\) is the indicator of a matching
\(M_0\), and
\(|M_0|=\sum_e p_e^2=\nu\).
Hence \(M_0\) is a maximum matching.

For every \(uv\in M_0\), the equality \(p_{uv}=1\)
implies \(L(u)=L(v)\).
For every edge outside \(M_0\), the corresponding lists
are disjoint. Admissibility gives disjointness for
nonadjacent vertices as well. Therefore the common lists
on the matched pairs and the lists on the unmatched
vertices are pairwise disjoint across distinct parts.

Conversely, any assignment of this form is admissible.
Each matched pair contributes \(\binom m2\) choices for
its color set, and each unmatched vertex contributes
\(m\) choices. Since different parts have disjoint lists,
the total number of color sets is
\(\binom m2^\nu m^{n-2\nu}\).
This proves the asserted characterization of the
minimizing assignments and completes the proof for
\(n\ge2\).

If \(n=1\), every color in the unique list gives a distinct
color set, so \(f'(G,m)=m\).
\end{proof}

\section{Odd-cycle blow-ups}\label{sec:blowups}

Fix an odd integer \(\ell\ge5\). Throughout this section,
\(t\ge1\) and \(m\ge\ell t\) are integers.
\subsection{Uniform counting estimates}

Using the entropy function \(h\) from
Section~\ref{sec:preliminaries}, put
\(q_0=\frac{\ell m}{2}, \qquad q=\lceil q_0\rceil, \qquad A=q_0h(2t/m)\).
We first record the uniform comparison
\begin{equation}\label{eq:uniform-stirling}
 \binom q{\ell t}
 =\Theta_\ell\!\left(t^{-1/2}e^A\right).
\end{equation}
Indeed, since
\(0<\frac{\ell t}{q}\le\frac{2t}{m}\le\frac2\ell<1\),
Lemma~\ref{lem:stirling} gives, uniformly in \(t\) and \(m\),
\[
 \binom q{\ell t}
 =\Theta_\ell\!\left(
 t^{-1/2}\exp\bigl(qh(\ell t/q)\bigr)\right).
\]
For completeness, write \(n=\ell t\) and keep
the prefactor from Stirling's formula explicitly.
It satisfies
\[
 \frac{1}{\sqrt{2\pi\ell t}}
 \le
 \sqrt{\frac{q}{2\pi n(q-n)}}
 =\frac{1}{\sqrt{2\pi n(1-n/q)}}
 \le
 \frac{1}{\sqrt{2\pi(\ell-2)t}}.
\]
Here \(1-2/\ell\le1-n/q<1\), so both bounding
constants are independent of the ratio \(m/t\).
The remaining factor in the exact Stirling formula
is \(\exp(\rho_q-\rho_n-\rho_{q-n})\).
Lemma~\ref{lem:stirling} bounds this factor above
and below by positive absolute constants, since
all three factorial arguments are positive integers.
Thus the prefactor is uniformly of order \(t^{-1/2}\),
including when the list size divided by the class size
tends to infinity.

To compare the exponents, define
\(\Phi(x)=xh(\ell t/x)\) for \(x>\ell t\).
For \(q_0\le x\le q\), we have
\(\Phi'(x) =-\log(1-\ell t/x) \le\log\frac{\ell}{\ell-2}\).
Since \(\Phi(q_0)=A\) and \(0\le q-q_0\le1/2\), it follows that
\[
 0\le qh(\ell t/q)-A
 \le\frac12\log\frac{\ell}{\ell-2}.
\]
This proves the comparison.
Consequently, Theorem~\ref{thm:blowup} is equivalent to
\[
 f'(G_{\ell,t},m)
 =\Theta_\ell\!\left(t^{-1/2}e^A\right).
\]
The upper bound follows from Lemma~\ref{lem:palettes}.
For the lower bound, we first control how many proper
colorings can use one prescribed set of colors, uniformly
over the admissible assignments. This is the multiplicity
estimate that will be paired with the injective-choice
bound from Lemma~\ref{lem:injective}.

\begin{lemma}\label{lem:fiber}
There exists a constant \(C_\ell>0\), depending only on
\(\ell\), such that for every
\(L\in\mathcal L_m(G_{\ell,t})\) and every set \(S\) of colors,
the number of proper \(L\)-colorings with image \(S\) is at most
\(C_\ell\sqrt t\,(2t/e)^{\ell t}\).
\end{lemma}

\begin{proof}
Let \(M_L(S)\) denote the number of proper \(L\)-colorings
with image \(S\). We may assume that \(M_L(S)>0\).
Since \(L\) is admissible, every proper \(L\)-coloring
is injective, and hence \(|S|=\ell t\).

For each \(c\in S\), its support
\(U_c=\{v\in V(G_{\ell,t}):c\in L(v)\}\)
is a clique. Every clique in \(G_{\ell,t}\) is contained
in one of the cliques
\(Q_i=V_i\cup V_{i+1}, \qquad i\in\mathbb Z/\ell\mathbb Z\).
Choose such a clique \(Q_i\) for each color \(c\in S\),
and make \(c\) available at every vertex of \(Q_i\).
Denote the resulting auxiliary list assignment by \(L^+\);
the occurrences of colors outside \(S\) are left unchanged.
Thus each color in the prescribed image has support exactly
the clique chosen for it. This enlargement can only increase
\(M_L(S)\), since each original coloring with that image
is still allowed, and hence \(M_L(S)\le M_{L^+}(S)\).
The lists in the auxiliary assignment may have more than
\(m\) elements, so \(L^+\) need not belong to
\(\mathcal L_m(G_{\ell,t})\).
This causes no difficulty: we use the enlargement only
to count colorings with the fixed image \(S\), and do
not apply any estimate requiring equal list sizes to it.
Its color supports remain cliques, so every proper
coloring of the auxiliary assignment is still injective.

Let \(x_i\) be the number of colors assigned to \(Q_i\).
Then
\(x_i\ge0, \qquad \sum_{i=0}^{\ell-1}x_i=\ell t\).
Independently for each color assigned to \(Q_i\), send it
to \(V_i\) or \(V_{i+1}\), each with probability \(1/2\).
Let \(Y_i\) count the colors sent from \(Q_i\) to \(V_i\).
The variables \(Y_i\) are independent, with
\(Y_i\sim\operatorname{Bin}(x_i,1/2)\).
The number of colors reaching \(V_i\) is
\(Z_i=Y_i+x_{i-1}-Y_{i-1}\).
Thus an allocation sends exactly \(t\) colors to every
class precisely on the event
\(E=\{Z_i=t\text{ for every }i\in\mathbb Z/\ell\mathbb Z\}\).
There are \(2^{\ell t}\) equally likely allocations.
Each successful allocation gives \((t!)^\ell\)
bijections from its assigned colors to the vertices
of the corresponding classes.
Conversely, a coloring of the auxiliary assignment
with image \(S\) determines the receiving class of
each color, and then the bijection within each class,
uniquely. Consequently, this allocation count is exact
for the auxiliary assignment:
\[
 M_{L^+}(S)=2^{\ell t}(t!)^\ell\PP(E).
\]
The monotonicity established above therefore gives
\[
 M_L(S)\le2^{\ell t}(t!)^\ell\PP(E).
\]
We prove that
\(\PP(E)=O_\ell\!\left(t^{-(\ell-1)/2}\right)\)
uniformly over all choices of \(x_0,\ldots,x_{\ell-1}\)
with sum \(\ell t\).

Suppose first that \(x_i\ge t/2\) for every \(i\).
Conditional on \(Y_0=y\), the equations
\(Y_i+x_{i-1}-Y_{i-1}=t, \qquad i=1,\ldots,\ell-1\),
determine \(Y_1,\ldots,Y_{\ell-1}\) successively.
By independence,
\[
 \PP(E\mid Y_0=y)
 \le
 \prod_{i=1}^{\ell-1}
       \max_{a\in\mathbb Z}\PP(Y_i=a).
\]
The binomial atom bound in Lemma~\ref{lem:stirling}
gives an absolute constant \(C>0\) such that
\(\max_{a\in\mathbb Z} \PP\bigl(\operatorname{Bin}(x,1/2)=a\bigr) \le Cx^{-1/2} \qquad (x\ge1)\).
Since each \(x_i\ge t/2\), averaging over \(Y_0\) yields
\(\PP(E)=O_\ell\!\left(t^{-(\ell-1)/2}\right)\).

Now suppose that \(x_j<t/2\) for some \(j\).
Put
\(\delta_i=x_{i-1}+x_i-2t\).
Since \(\ell\) is odd, cancellation in the alternating
sum gives
\[
 \sum_{r=0}^{\ell-1}(-1)^r\delta_{j-r}
 =2(x_j-t).
\]
The absolute value of the right-hand side is greater
than \(t\). Hence there is an index \(i\) for which
\(|\delta_i|>\frac{t}{\ell}\).
Put \(N_i=x_{i-1}+x_i\). Since
\(x_{i-1}-Y_{i-1}\) has distribution
\(\operatorname{Bin}(x_{i-1},1/2)\) and is independent
of \(Y_i\), we have
\(Z_i\sim\operatorname{Bin}(N_i,1/2)\).
Moreover,
\[
 \left|t-\frac{N_i}{2}\right|
 =\frac{|\delta_i|}{2}
 >\frac{t}{2\ell},
 \qquad
 N_i\le\ell t.
\]
If \(N_i=0\), then \(Z_i=0\) almost surely and
\(\PP(E)=0\), since \(t\ge1\).
Otherwise, Lemma~\ref{lem:hoeffding} gives
\[
 \begin{aligned}
 \PP(E)
 &\le\PP(Z_i=t)\\
 &\le
 \PP\left(\left|Z_i-\frac{N_i}{2}\right|
                   \ge\frac{t}{2\ell}\right)\\
 &\le2\exp\left(-\frac{t^2}{2\ell^2N_i}\right)\\
 &\le2\exp\left(-\frac{t}{2\ell^3}\right).
 \end{aligned}
\]
For fixed \(\ell\), this is
\(O_\ell(t^{-(\ell-1)/2})\) for all \(t\ge1\).
Thus the required probability bound holds in both cases.

Finally, Lemma~\ref{lem:stirling} gives
\(t!\le C\sqrt t\,(t/e)^t \qquad (t\ge1)\)
for an absolute constant \(C\).
Combining the preceding estimates, we obtain
\[
 \begin{aligned}
 M_L(S)
 &\le2^{\ell t}(t!)^\ell\PP(E)\\
 &\le
 C_\ell\,2^{\ell t}t^{\ell/2}(t/e)^{\ell t}
                         t^{-(\ell-1)/2}\\
 &=C_\ell\sqrt t\,(2t/e)^{\ell t},
 \end{aligned}
\]
after choosing \(C_\ell\) sufficiently large.
\end{proof}

The preceding lemma bounds the count for one image,
while Lemma~\ref{lem:injective} estimates the total number
of injective choices. Their comparison cancels the linear
terms in the exponents, leaving the entropy term, the
square-root factor, and a nonnegative support-dependent
remainder. We retain that remainder for the stability
argument.

\begin{proof}[Proof of Theorem~\ref{thm:blowup}]
Put
\(n=\ell t,\qquad r=2t,\qquad p=\frac rm, \qquad q=\left\lceil\frac{\ell m}{2}\right\rceil\).
Then \(0<p\le2/\ell<1\), and the exponent defined at the
beginning of this section satisfies
\(A=\frac np h(p) =\frac{\ell m}{2}h(2t/m)\).
By the uniform Stirling comparison established above,
\(\binom qn =\Theta_\ell\!\left(t^{-1/2}e^A\right)\).
Lemma~\ref{lem:palettes} gives
\(f'(G_{\ell,t},m)\le\binom qn\).
It therefore remains to prove a lower bound of the form
\(c_\ell t^{-1/2}e^A\), where \(c_\ell>0\) depends only
on \(\ell\).
Fix \(L\in\mathcal L_m(G_{\ell,t})\), and let
\(\mathcal U=\bigcup_{v\in V(G_{\ell,t})}L(v)\).
For each \(c\in\mathcal U\), write
\(d_c=|\{v:c\in L(v)\}|, \qquad p_c=\frac{d_c}{m}\).
The support of each color is a clique, so
\(1\le d_c\le\omega(G_{\ell,t})=2t=r<m\).
Consequently,
\(0<p_c\le p<1, \qquad \sum_{c\in\mathcal U}p_c =\frac1m\sum_{c\in\mathcal U}d_c=n\).

Let \(T_L\) be the number of proper \(L\)-colorings, and
put \(N_L=|\mathcal C_{G_{\ell,t}}(L)|\).
Recall the function
\(g(z)=-\frac{(1-z)\log(1-z)}{z}, \qquad 0<z<1\),
which is decreasing, as observed in the proof of
Proposition~\ref{prop:general}.
Lemma~\ref{lem:injective} applies because \(d_c<m\)
for every \(c\), and gives
\[
 \log T_L
 \ge n\log m-n
       +\sum_{c\in\mathcal U}p_cg(p_c).
\]
On the other hand, Lemma~\ref{lem:fiber} bounds the number
of colorings with any prescribed image by
\(C_\ell\sqrt t\,(r/e)^n\).
Summing over the \(N_L\) possible images, we obtain
\(T_L\le N_L C_\ell\sqrt t\,(r/e)^n\).
Combining these inequalities and using
\(-\log p+g(p)=\frac{h(p)}p\),
we obtain
\begin{equation}\label{eq:slack}
 \begin{aligned}
 \log N_L
 &\ge
 n\log(m/r)
 +\sum_{c\in\mathcal U}p_cg(p_c)
 -\frac12\log t-\log C_\ell\\
 &=
 \frac np h(p)-\frac12\log t-\log C_\ell
 +\sum_{c\in\mathcal U}
       p_c\bigl(g(p_c)-g(p)\bigr).
 \end{aligned}
\end{equation}
Here the terms \(-n\) in the lower bound for
\(\log T_L\) and in the logarithm of the multiplicity
bound cancel.

Since \(p_c\le p\) and \(g\) is decreasing, every summand
in the last sum of \eqref{eq:slack} is nonnegative.
Thus
\(N_L\ge C_\ell^{-1}t^{-1/2}e^A\).
As \(L\) was arbitrary, minimizing over
\(\mathcal L_m(G_{\ell,t})\) yields
\[
 f'(G_{\ell,t},m)
 \ge C_\ell^{-1}t^{-1/2}e^A.
\]
Together with the upper bound and the uniform Stirling
comparison, this proves
\[
 f'(G_{\ell,t},m)
 =\Theta_\ell\!\left(
      \binom{\lceil\ell m/2\rceil}{\ell t}\right).
\]
All constants depend only on \(\ell\), so the estimate
is uniform over \(t\ge1\) and \(m\ge\ell t\).
\end{proof}

\subsection{The gap from clique partitions}

Comparing the uniform counting estimate with the exact
clique-partition value gives both the exponential term
and the polynomial factor in their ratio. The resulting
bound remains uniform over the allowed list sizes and
quantifies how far clique partitions are from optimal.

\begin{corollary}\label{cor:uniform-gap}
Fix an odd integer \(\ell\ge5\). Uniformly over integers
\(t\ge1\) and \(m\ge\ell t\),
\[
 \log\frac{w(G_{\ell,t},m)}{f'(G_{\ell,t},m)}
 =
 m\left(h(t/m)-\frac12h(2t/m)\right)
 -\frac{\ell-1}{4}\log t+O_\ell(1).
\]
In particular, there is a constant \(c_\ell>0\) such that
\[
 \frac{w(G_{\ell,t},m)}{f'(G_{\ell,t},m)}
 \ge c_\ell\,2^t t^{-(\ell-1)/4}
\]
for all such \(t\) and \(m\).
Consequently, there exists \(t_0=t_0(\ell)\) such that
\[
 f'(G_{\ell,t},m)<w(G_{\ell,t},m)
 \qquad
 (t\ge t_0,\ m\ge\ell t).
\]
In particular, \(f(G_{\ell,t})<w(G_{\ell,t})\) for
every \(t\ge t_0\).
\end{corollary}

\begin{proof}
Write \(\ell=2k+1\).
By Lemma~\ref{lem:palettes},
\(w(G_{\ell,t},m)=\binom m{2t}^{k}\binom mt\).
Since \(2t/m\le2/\ell<1\), Lemma~\ref{lem:stirling} gives,
uniformly for \(m\ge\ell t\),
\[
 \log w(G_{\ell,t},m)
 =
 kmh(2t/m)+mh(t/m)
 -\frac{k+1}{2}\log t+O_\ell(1).
\]
On the other hand, Theorem~\ref{thm:blowup} and the
uniform Stirling comparison give
\[
 \log f'(G_{\ell,t},m)
 =
 \frac{\ell m}{2}h(2t/m)
 -\frac12\log t+O_\ell(1).
\]
Subtracting these expressions and using
\(k=(\ell-1)/2\) proves the first assertion.

Put \(x=t/m\), so that \(0<x\le1/\ell\), and define
\(u(z)=-(1-z)\log(1-z), \qquad 0\le z<1\).
The function \(u\) is concave and satisfies \(u(0)=0\).
Hence
\(u(x)\ge\frac12u(2x)\).
It follows that
\[
 \begin{aligned}
 h(x)-\frac12h(2x)
 &=x\log2+u(x)-\frac12u(2x)\\
 &\ge x\log2.
 \end{aligned}
\]
The first assertion therefore yields a constant
\(C_\ell\) such that
\[
 \log\frac{w(G_{\ell,t},m)}{f'(G_{\ell,t},m)}
 \ge
 t\log2-\frac{\ell-1}{4}\log t-C_\ell.
\]
Exponentiating proves the stated lower bound with
\(c_\ell=e^{-C_\ell}\).

Since \(2^t t^{-(\ell-1)/4}\) tends to infinity,
the ratio is greater than one for all sufficiently
large \(t\), with a threshold depending only on
\(\ell\). Taking \(m=\ell t\) and applying
Lemma~\ref{lem:join} gives the final assertion.
\end{proof}

For the five-cycle blow-ups, take \(m=5t\).
Since \(f(G_{5,t})=f'(G_{5,t},5t)\) by
Lemma~\ref{lem:join}, Theorem~\ref{thm:blowup} and
Corollary~\ref{cor:uniform-gap} give
\[
 \begin{aligned}
 \lim_{t\to\infty}\frac{\log f(G_{5,t})}{5t}
 &=\frac52h(2/5),\\
 \lim_{t\to\infty}\frac1{5t}
       \log\frac{w(G_{5,t})}{f(G_{5,t})}
 &=h(1/5)-\frac12h(2/5).
 \end{aligned}
\]
The second limit is positive; indeed, the entropy
inequality in the preceding proof shows that it is
at least \((\log2)/5\).

\subsection{Stability of near-minimizing assignments}

The support-dependent remainder penalizes colors whose
supports are smaller than a maximal clique. The argument
below turns this penalty into an incidence edit bound:
large supports are first completed to consecutive pairs
of classes, and the resulting palettes are then adjusted
to the allowed sizes prescribed in Lemma~\ref{lem:palettes}.
For list assignments \(L\) and \(\widetilde L\) on
the same graph, define their incidence distance by
\[
 d_{\mathrm{inc}}(L,\widetilde L)
 =
 \sum_{v\in V(G)}
       |L(v)\triangle\widetilde L(v)|.
\]
This is the minimum number of insertions and deletions
of vertex--color incidences needed to transform one
assignment into the other. New colors may be introduced
in \(\widetilde L\).

\begin{corollary}\label{cor:stability}
Fix an odd integer \(\ell\ge5\) and constants
\(K\ge1\) and \(M\ge\ell\).
Suppose that \(t,m\) are integers satisfying
\(\ell t\le m\le Mt\), and that
\(L\in\mathcal L_m(G_{\ell,t})\) satisfies
\[
 N_L:=|\mathcal C_{G_{\ell,t}}(L)|
 \le K\binom{\lceil\ell m/2\rceil}{\ell t}.
\]
Then there is an assignment \(\widetilde L\) of the
form constructed in Lemma~\ref{lem:palettes} such that
\[
 d_{\mathrm{inc}}(L,\widetilde L)
 \le C_\ell M^2(1+\log K)t,
\]
where \(C_\ell>0\) depends only on \(\ell\).
\end{corollary}

\begin{proof}
Put \(r=2t\) and \(p=r/m\).
Let \(\mathcal U\) be the union of the lists of \(L\),
and, for \(c\in\mathcal U\), put
\(d_c=|\{v:c\in L(v)\}|, \qquad p_c=d_c/m\).
Define
\[
 \Delta(L)
 =
 \sum_{c\in\mathcal U}
 p_c\bigl(g(p_c)-g(p)\bigr),
\]
where
\(g(z)=-\frac{(1-z)\log(1-z)}{z}, \qquad 0<z<1\).
We first specify the constants in the bound for
this remainder. Choose \(C_\ell^{\mathrm{fib}}\ge1\)
so that the multiplicity bound in Lemma~\ref{lem:fiber}
holds with that constant; then \eqref{eq:slack} holds
with the same constant in place of the unspecified one.
Choose also \(C_\ell^{\mathrm{St}}\ge1\), depending
only on \(\ell\), such that the uniform Stirling
comparison gives
\[
 \binom{\lceil\ell m/2\rceil}{\ell t}
 \le C_\ell^{\mathrm{St}}t^{-1/2}e^A.
\]
In terms of these constants, \eqref{eq:slack}
and the hypothesis on \(N_L\) give
\[
 \begin{aligned}
 A-\frac12\log t-\log C_\ell^{\mathrm{fib}}
       +\Delta(L)
 &\le\log N_L\\
 &\le A-\frac12\log t+\log K
                 +\log C_\ell^{\mathrm{St}}.
 \end{aligned}
\]
Every summand defining the remainder is nonnegative,
because the function being compared is decreasing.
It follows that
\[
 0\le\Delta(L)
 \le\log K+
       \log\bigl(C_\ell^{\mathrm{fib}}
                 C_\ell^{\mathrm{St}}\bigr).
\]
Define
\(B_\ell=\max\{1,\log(C_\ell^{\mathrm{fib}}C_\ell^{\mathrm{St}})\}\).
Then \(B_\ell\ge1\), it depends only on \(\ell\),
and we have the bound
\(0\le\Delta(L)\le B_\ell+\log K\).

We first convert this bound into control of the color
supports. For \(0<z<1\),
\[
 -g'(z)
 =
 \frac{-\log(1-z)-z}{z^2}
 =
 \frac1{z^2}\int_0^z\frac{u}{1-u}\,du
 \ge\frac12.
\]
Since \(0<p_c\le p<1\), it follows that
\(g(p_c)-g(p)\ge\frac12(p-p_c)\).
Consequently,
\[
 \Delta(L)
 \ge\frac12\sum_{c\in\mathcal U}p_c(p-p_c)
 =\frac1{2m^2}
       \sum_{c\in\mathcal U}d_c(2t-d_c),
\]
and hence
\(\sum_{c\in\mathcal U}d_c(2t-d_c) \le2m^2(B_\ell+\log K)\).

Delete every color with \(d_c\le t\).
For each remaining color, its support has more than
\(t\) vertices and is a clique. It therefore meets
exactly two consecutive classes and is contained in
a unique clique
\(Q_i=V_i\cup V_{i+1}\).
Enlarge its support to all of \(Q_i\).
Let \(L_0\) be the resulting list assignment, and put
\(D_0=d_{\mathrm{inc}}(L,L_0)\).
Its list sizes need not equal the original list size:
the incidence distance is defined for arbitrary finite
list assignments on the same graph, and both deletions
and insertions are included in it.
The intermediate assignment is used only to organize
the palettes and count these edits. Equal list sizes
will be restored in the final palette assignment.
The number of edits is
\(D_0 = \sum_{d_c\le t}d_c +\sum_{d_c>t}(2t-d_c)\).
For \(d_c\le t\), we have
\(d_c\le d_c(2t-d_c)/t\).
For \(d_c>t\), we likewise have
\(2t-d_c\le d_c(2t-d_c)/t\).
Thus
\[
 \begin{aligned}
 D_0
 &\le\frac1t
       \sum_{c\in\mathcal U}d_c(2t-d_c)\\
 &\le\frac{2m^2}{t}(B_\ell+\log K)\\
 &\le2M^2(B_\ell+\log K)t.
 \end{aligned}
\]

Let \(P_i\) be the palette of colors whose support in
\(L_0\) is \(Q_i\), and put \(x_i=|P_i|\).
These palettes are pairwise disjoint, and every
vertex in \(V_i\) has list \(P_{i-1}\cup P_i\).
Write
\(a_i=x_{i-1}+x_i\).
Since every original list has size \(m\),
\[
 \begin{aligned}
 t\sum_{i=0}^{\ell-1}|a_i-m|
 &=
 \sum_{v\in V(G_{\ell,t})}
       \bigl||L_0(v)|-|L(v)|\bigr|\\
 &\le d_{\mathrm{inc}}(L,L_0)
 =D_0.
 \end{aligned}
\]
The oddness of \(\ell\) gives the identity
\[
 2x_i-m
 =
 \sum_{j=0}^{\ell-1}(-1)^j(a_{i-j}-m),
\]
where indices are taken modulo \(\ell\).
Therefore,
\[
 \left|x_i-\frac m2\right|
 \le\frac12\sum_{j=0}^{\ell-1}|a_j-m|
 \le\frac{D_0}{2t}
 \qquad\text{for every }i.
\]

Let \(y_0,\ldots,y_{\ell-1}\) be the palette sizes
prescribed in Lemma~\ref{lem:palettes}, with exceptional
class \(V_{\ell-1}\) when \(m\) is odd.
In either parity,
\(\left|y_i-\frac m2\right|\le\frac12\).
Adjust each \(P_i\) to have size \(y_i\), deleting
colors or adding fresh colors as needed.
The palettes remain pairwise disjoint.
Each palette color has support \(Q_i\), so these
adjustments require at most
\[
 \begin{aligned}
 2t\sum_{i=0}^{\ell-1}|x_i-y_i|
 &\le
 2t\sum_{i=0}^{\ell-1}
 \left(
 \left|x_i-\frac m2\right|
 +\left|y_i-\frac m2\right|
 \right)\\
 &\le\ell D_0+\ell t
 \end{aligned}
\]
incidence edits.
If \(m\) is odd, add one further fresh color to the
list of every vertex in \(V_{\ell-1}\), using \(t\)
additional edits. The resulting assignment
\(\widetilde L\) is precisely of the form constructed
in Lemma~\ref{lem:palettes}. Its distance from \(L\)
satisfies
\[
 \begin{aligned}
 d_{\mathrm{inc}}(L,\widetilde L)
 &\le D_0+\ell D_0+\ell t+t\\
 &=(\ell+1)D_0+(\ell+1)t\\
 &\le C_\ell M^2(1+\log K)t,
 \end{aligned}
\]
where \(C_\ell\) may be chosen explicitly as
\(C_\ell=(\ell+1)(2B_\ell+1)\).
Indeed, the bound on the intermediate distance gives
\[
 \begin{aligned}
 (\ell+1)D_0+(\ell+1)t
 &\le(\ell+1)
       \bigl(2M^2(B_\ell+\log K)+1\bigr)t\\
 &\le(\ell+1)M^2
       \bigl(2B_\ell+2\log K+1\bigr)t\\
 &\le(\ell+1)(2B_\ell+1)M^2(1+\log K)t.
 \end{aligned}
\]
The second inequality uses \(M^2\ge1\), and the last
uses \(B_\ell\ge1\) and \(\log K\ge0\).
Thus the constant in the stated distance estimate
depends only on the fixed cycle length, as required.
\end{proof}

\section{Equality in the square lower bound}\label{sec:star}

We finish by considering equality in the square lower
bound. Theorem~\ref{thm:known-bounds} records the bound
and the examples given by Gravier, Maffray, and
Mohar~\cite[Theorem~1 and Claim~1]{GMM}.
Sufficiency is therefore already known. The task here
is to prove the converse, showing that every noncomplete
graph attaining the bound belongs to that same family.

\begin{theorem}\label{thm:star}
Let \(G\) be a noncomplete graph, and put
\(n=|V(G)|\ge2\).
Then \(f(G)=n^2\) if and only if all edges of
\(\overline G\) have a common endpoint.
\end{theorem}

\begin{proof}
The lower bound and the sufficiency of the stated condition
are the cited assertions in Theorem~\ref{thm:known-bounds}.
We prove necessity.
The case \(n=2\) is immediate, so assume \(n\ge3\).
Suppose that \(f(G)=n^2\).
By Lemma~\ref{lem:join}, there is an assignment
\(L\in\mathcal L_n(G)\) such that
\(N_L:=|\mathcal C_G(L)|=n^2\).
Let \(u,v\) be any nonadjacent vertices, and put
\(U=L(u)\cup L(v)\).
Admissibility gives \(L(u)\cap L(v)=\varnothing\),
so \(|U|=2n\).

We first show that every list is contained in \(U\).
Suppose instead that some
\(z\in V(G)\setminus\{u,v\}\) has a color
\(c\in L(z)\setminus U\).
On \(H=G-z\), define
\(J(w)=L(w)\setminus\{c\} \qquad (w\in V(H))\).
Every list of \(J\) has size at least \(n-1=|V(H)|\),
while
\(J(u)=L(u),\qquad J(v)=L(v)\).
Lemma~\ref{lem:two-list-bound} therefore gives at least
\(n^2\) distinct color sets of proper \(J\)-colorings of
\(H\). None contains \(c\).
Extending each such coloring by assigning \(c\) to \(z\)
produces at least \(n^2\) distinct color sets of proper
\(L\)-colorings containing \(c\).

There is also a proper \(L\)-coloring avoiding \(c\).
To see this, put
\(L^-(w)=L(w)\setminus\{c\} \qquad (w\in V(G))\).
For every nonempty proper subset \(W\subsetneq V(G)\),
\(\left|\bigcup_{w\in W}L^-(w)\right| \ge n-1\ge |W|\).
For \(W=V(G)\), the union contains the unchanged,
disjoint lists \(L(u)\) and \(L(v)\), and hence has
size at least \(2n\).
Theorem~\ref{thm:hall} therefore gives a system of distinct
representatives of the lists \(L^-(w)\).
This is an injective proper \(L\)-coloring avoiding \(c\).
Together with the preceding \(n^2\) sets containing \(c\),
it contradicts \(N_L=n^2\).

Thus, for every nonedge \(uv\),
\[
 \bigcup_{w\in V(G)}L(w)=L(u)\cup L(v),
 \qquad
 |L(u)\cup L(v)|=2n.
\]
In particular, \(\overline G\) is triangle-free:
a triangle in \(\overline G\) would give three pairwise
disjoint \(n\)-element lists contained in a set of
size \(2n\).

Suppose now that \(\overline G\) contains two disjoint
edges \(uv\) and \(xy\).
Write
\(A=L(u),\qquad B=L(v),\qquad U=A\cup B\).
The lists \(L(x)\) and \(L(y)\) are disjoint
\(n\)-element subsets of \(U\), so
\(L(x)\cup L(y)=U\).
Choose \(c\in A\), and let \(z\in\{x,y\}\) be the
unique vertex whose list contains \(c\).
On \(G-z\), delete \(c\) from every remaining list.
The resulting lists have size at least \(n-1\),
and the lists at \(u,v\) are the disjoint sets
\(A\setminus\{c\}\) and \(B\), of sizes \(n-1\)
and \(n\), respectively.
Lemma~\ref{lem:two-list-bound}, followed by extension
with color \(c\) at \(z\), gives at least
\(n(n-1)\)
distinct color sets of proper \(L\)-colorings
containing \(c\).

We next count color sets avoiding \(c\).
Since \(L(x)\cup L(y)=U\), at least one vertex
\(w\in\{x,y\}\) has a color \(d\in L(w)\cap B\).
Assign \(d\) to \(w\), and then color the remaining
vertices of \(G-\{u,v\}\) greedily, always avoiding
\(c\) and all previously chosen colors.
This is possible because each list after deleting
\(c\) has size at least \(n-1\), whereas the subgraph
being colored has only \(n-2\) vertices.
Let \(T\) be the image of this fixed injective coloring.
Then
\(|T|=n-2,\qquad T\subseteq U\setminus\{c\},\qquad T\cap B\ne\varnothing\).
Put
\(a=|T\cap A|,\qquad b=|T\cap B|\).
Thus \(a+b=n-2\) and \(b\ge1\).
The available colors at \(u,v\), while still avoiding
\(c\), form the disjoint sets
\(A\setminus(T\cup\{c\}), \qquad B\setminus T\).
Their sizes are
\(\lambda_u=n-1-a,\qquad \lambda_v=n-b\).
Consequently,
\(2\le\lambda_u,\lambda_v\le n-1, \qquad \lambda_u+\lambda_v=n+1\).
It follows that
\[
 \begin{aligned}
 \lambda_u\lambda_v
 &=
 2(\lambda_u+\lambda_v)-4
 +(\lambda_u-2)(\lambda_v-2)\\
 &\ge2(n-1).
 \end{aligned}
\]

Each pair of available colors extends the fixed coloring
to a proper \(L\)-coloring avoiding \(c\).
The resulting color sets are distinct: both new colors
lie outside \(T\), and their membership in the disjoint
sets \(A\) and \(B\) determines which was assigned
to \(u\) and which to \(v\).
There are therefore at least \(2(n-1)\) color sets
avoiding \(c\).
The sets containing \(c\) and those avoiding \(c\)
are disjoint, so
\[
 N_L
 \ge n(n-1)+2(n-1)
 =n^2+n-2>n^2,
\]
where \(n\ge4\) because the two nonedges are disjoint.
This contradiction shows that the edges of
\(\overline G\) pairwise intersect.

Finally, \(\overline G\) is nonempty and triangle-free.
If it has only one edge, that edge has a common endpoint
with every edge. Otherwise, choose two distinct edges
\(ab\) and \(ac\).
Any edge avoiding \(a\) must meet both of these edges,
and hence must be \(bc\), which would create a triangle.
Thus every edge contains \(a\), proving the required
common-endpoint condition.
\end{proof}

\section*{Declaration of generative AI and AI-assisted technologies}
The author used a large language model during the preparation of this manuscript for assistance with mathematical exploration, checking arguments, organization and language editing. The author is responsible for the mathematical statements, proofs, references and final text.

\appendix

\section{Counting color sets by exact supports}
\label{app:support-count}

For a fixed admissible assignment, Hall's condition
depends on a color only through its support. Grouping
the colors by support separates feasibility from the
choice of individual colors. The resulting formula
will be used in the finite evaluation in
Appendix~\ref{sec:finite-count}. By Theorem~\ref{thm:hall},
the feasible profiles are determined entirely by the
support incidences; the binomial weights then count
all the color sets represented by each such profile.

\begin{proposition}[A support-count form of Hall's theorem]\label{prop:support-count}
Let \(G\) be a graph, put \(n=|V(G)|\), and let
\(m\ge n\) and \(L\in\mathcal L_m(G)\).
Let \(\mathcal K(G)\) be the family of nonempty cliques
of \(G\). For each color occurring in \(L\), define
\(\operatorname{supp}_L(c) =\{v\in V(G):c\in L(v)\}\),
and, for \(Q\in\mathcal K(G)\), put
\(a_Q= |\{c:\operatorname{supp}_L(c)=Q\}|\).
Then
\(\sum_{\substack{Q\in\mathcal K(G)\\v\in Q}}a_Q=m
\qquad (v\in V(G))\).
Moreover,
\[
 |\mathcal C_G(L)|
 =
 \sum_{x\in\mathcal F(a)}
       \prod_{Q\in\mathcal K(G)}\binom{a_Q}{x_Q},
\]
where \(\mathcal F(a)\) consists of the integer vectors
\(x=(x_Q)_{Q\in\mathcal K(G)}\) satisfying
\[
 \begin{gathered}
 0\le x_Q\le a_Q
       \quad(Q\in\mathcal K(G)),
 \qquad
 \sum_{Q\in\mathcal K(G)}x_Q=n,\\
 \sum_{\substack{Q\in\mathcal K(G)\\
                  Q\cap U\ne\varnothing}}x_Q
 \ge |U|
       \qquad(U\subseteq V(G)).
 \end{gathered}
\]
\end{proposition}

\begin{proof}
Admissibility implies that the support of every color
is a nonempty clique. Grouping the colors according
to their supports therefore gives the stated identity
for each list size.
Consider an \(n\)-element set \(S\) of colors, and let
\(x_Q\) count its colors with support \(Q\).
For \(U\subseteq V(G)\), its neighborhood in the
incidence graph between \(V(G)\) and \(S\) has size
\(\sum_{\substack{Q\in\mathcal K(G)\\ Q\cap U\ne\varnothing}}x_Q\).
Theorem~\ref{thm:hall} thus shows that \(S\) is realizable
if and only if its profile belongs to
\(\mathcal F(a)\).

For each feasible profile \(x\), independently choosing
\(x_Q\) colors from the \(a_Q\) colors supported on
\(Q\) gives
\(\prod_{Q\in\mathcal K(G)}\binom{a_Q}{x_Q}\)
distinct color sets. Every color set has a unique
profile, so summing over the feasible profiles proves
the formula.
\end{proof}

\section{An exact count for the 20-vertex construction}
\label{sec:finite-count}

Proposition~\ref{prop:support-count} expresses the exact
count through feasible support profiles and their binomial
weights. In the construction of Lemma~\ref{lem:palettes},
the identical lists within each class reduce Hall's
conditions to a small family of inequalities in the five
palette counts. Retaining only the profiles of realizable
sets gives a sharper upper bound for \(G_{5,4}\), and the
resulting finite sum can be evaluated exactly. We first
specify the palettes and lists used in this calculation.
Let \(P_0,\ldots,P_4\) be pairwise disjoint palettes,
each of size \(10\), and define
\(L^\ast(v)=P_{i-1}\cup P_i \qquad (v\in V_i)\),
where indices are taken modulo \(5\).
Write
\(\mathcal P=\bigcup_{i=0}^4P_i\).
Thus \(|\mathcal P|=50\), and every list has size \(20\).

For a \(20\)-element set \(S\subseteq\mathcal P\), put
\(x_i=|S\cap P_i|\).
Since \(L^\ast\) is admissible, every proper
\(L^\ast\)-coloring is injective.
Consequently, \(S\) is realizable if and only if the
bipartite incidence graph between \(V(G_{5,4})\) and
\(S\), with an edge \(vc\) whenever \(c\in L^\ast(v)\),
has a matching saturating all vertices. By
Theorem~\ref{thm:hall}, this is equivalent to Hall's
inequalities for its vertex subsets.
It suffices to check Hall's condition for unions of
whole classes. Indeed, if a vertex set meets \(V_i\),
adding the remaining vertices of \(V_i\) leaves its
color neighborhood unchanged, since all vertices in
\(V_i\) have the same list.
Thus, for \(J\subseteq\mathbb Z/5\mathbb Z\), the
required conditions are
\[
 \sum_{j\in J\cup(J-1)}x_j\ge4|J|,
 \qquad
 J-1=\{i-1:i\in J\}.
\]

A single class gives \(x_{i-1}+x_i\ge4\).
Two consecutive classes see three consecutive palettes.
Since \(\sum_i x_i=20\), their Hall condition is
equivalent to requiring at most \(12\) selected colors
in the complementary pair of consecutive palettes.
Three consecutive classes see all but one palette,
so their Hall conditions are precisely \(x_i\le8\)
for every \(i\).
These inequalities also imply all remaining Hall
conditions. Two nonconsecutive classes see four
palettes, containing at least \(20-8=12\) selected
colors, which is sufficient for their eight vertices.
Three classes that are not consecutive see all five
palettes, as does every union of four or five classes.
Their Hall conditions follow from \(|S|=20\).

Therefore the feasible profiles are exactly
\[
 \mathcal X=
 \left\{
 x\in\{0,\ldots,8\}^{5}:
 \sum_{i=0}^4x_i=20,\quad
 4\le x_{i-1}+x_i\le12
 \text{ for every }i
 \right\}.
\]
For a fixed profile \(x\), there are
\(\prod_{i=0}^4\binom{10}{x_i}\) choices of \(S\).
Hence
\[
 \begin{aligned}
 N_{L^\ast}
 &:=
 |\mathcal C_{G_{5,4}}(L^\ast)|\\
 &=\sum_{x\in\mathcal X}
       \prod_{i=0}^4\binom{10}{x_i}\\
 &=45\,699\,529\,371\,705.
 \end{aligned}
\]
For direct evaluation, one may enumerate
\(x_0,\ldots,x_3\in\{0,\ldots,8\}\), set
\(x_4=20-\sum_{i=0}^3x_i\), and retain the vectors
satisfying the displayed constraints.
This gives \(|\mathcal X|=2101\).
By Lemma~\ref{lem:join}, we obtain the improved bound
\[
 f(G_{5,4})
 =f'(G_{5,4},20)
 \le45\,699\,529\,371\,705
 <\binom{50}{20}.
\]

\end{document}